\documentstyle[11pt,twoside, amscd]{amsart}

\newcommand{\Hilb}{\operatorname{Hilb}}

\newcommand{\red}{\operatorname{red}}

\newcommand{\id}{\operatorname{id}}

\newtheorem{thm}{Theorem}[section]

\newtheorem{lem}[thm]{Lemma}
\newtheorem{cor}[thm]{Corollary}

\newtheorem{prop}[thm]{Proposition}

\theoremstyle{definition}
\newtheorem{defn}[thm]{Definition}

\newtheorem{say}[thm]{}

\newtheorem{rem}[thm]{Remark}           %\renewcommand{\theremark}{}

\theoremstyle{remark}

\begin{document}
\title[Relative GIT and Universal Moduli Spaces]{Relative Geometric
Invariant Theory and Universal Moduli Spaces}
\author[Yi Hu]{Yi Hu$^*$}
\date{April 16, 1995}
%\subjclass
%Primary 14D25, 14L30,  Secondary 14C20, 14A99 \endsubjclass
\thanks{$^*$Research partially supported by NSF grant DMS 9401695.}
\maketitle

{\footnotesize
%{\scriptsize
\tableofcontents
}

\section{Introduction}
\label{sec:introduction}
\begin{say} {\sl Motivation: the universal moduli problems}.
The motivation of this paper is to lay a GIT ground and to apply it to the
so-called universal
moduli problems such as the following ones.
\begin{enumerate}
\item The universal moduli space $\overline{FM_{g, n}} \rightarrow
\overline{M_g}$
      of Fulton-MacPherson configuration spaces of stable curves. That is,
given a stable
      curve $[C]$ in  $\overline{M_g}$, the fiber in the universal moduli space
will be
      $C[n]/\hbox{Aut}(C)$,
      the Fulton-MacPherson configuration space for the curve $C$ modulo the
automorphism
      group of $C$.  (See also \cite{Pandharipande94b})

\item The compactified
      universal Picard $\overline{P^d_g}  \rightarrow \overline{M_g}$ of degree
$d$ line bundles
      (\cite{Caporaso94}).

\item The universal moduli space $\overline{P_{g,m}(e, r, F, \alpha)}
\rightarrow
     \overline{M_{g,n}}$ of $p$-semistable parabolic sheaves of degree $e$,
rank $r$,
     type $F$, and weight $\alpha$  (\cite{Hu95}).

\item The universal moduli space $M_{g}({\cal O}, P) \rightarrow
\overline{M_g}$
of $p$-semistable coherent sheaves of pure dimension 1
 with a fixed Hilbert polynomial $P$ such that
the fiber over a stable curve $[C]$ is  functorially identified with
Simpson's moduli space $M_{C}({\cal O}_C, P)$ of $p$-semistable coherent
sheaves over $C$
of the Hilbert polynomial $P$ modulo the automorphism group of $C$.

The construction of the universal moduli space $M_{g}({\cal O}, P)$, as
described above,
bears a straightforward generalization to the universal moduli over
 the moduli spaces of higher dimensional varieties (e.g.,
  surfaces of general type, Clabi-Yau 3-folds, etc, see \cite{Hu95}) because
Simpson's construction
  \cite{Simpson94} works for any projective scheme.

\item The universal Hilbert scheme $\hbox{Hilb}^n_g  \rightarrow
\overline{M_g}$
of 0-dimensional subschemes of length $n$ on the Mumford-Deligne  stable curves
such that the fiber over a stable curve $[C]$ is canonically the Hilbert scheme
of 0-dimensional subschemes of length $n$ on $C$ modulo the automorphism group
of $C$.
Moreover, when $P(x) = x + n + 1 - g$ there exists a canonical dominating
morphism $$\psi: \hbox{Hilb}^n_g  \rightarrow M_{g}({\cal O}, P).$$
\end{enumerate}
\end{say}

\begin{say} {\sl Relative GIT}.
   All these universal moduli problems correspond to the following GIT setup
which we shall
call Relative GIT.

   We have a projective map $$\pi: Y \longrightarrow X$$
   and an epimorphism of two reductive algebraic groups
   $$\rho: G' \longrightarrow G.$$ We assume that $\pi: Y \longrightarrow X$
   is equivariant with respect to the  homomorphism $\rho: G' \longrightarrow
G$.

\medskip\noindent
{\bf  Question RGIT}.  Given a linearization $L$ on $X$ and the
        GIT quotient $X^{ss}(L)/\!/G$, find a  linearization $M$ on $Y$ and the
 GIT quotient
      $Y^{ss}(M)/\!/G'$ so that  $Y^{ss}(M)/\!/G'$ factors naturally to
$X^{ss}(L)/\!/G$.
\medskip

Using the relation between GIT and moment maps, we solved in this paper the
Question RGIT.
In the  case when $X^{ss} = X^s$ (which is the case for $\overline{M_g}$),
the question has a particularly nice solution: the stable locus $X^{ss}=X^s$
downstairs essentially
determines the stable locus $Y^{ss}(M) =Y^s(M)$ upstair. These results are
contained
in \S\S 3 and 4.

The solution to Question RGIT leads to a unified and easy approach to all
{\it those} universal moduli problems where the base moduli spaces
have satisfactory GIT constructions (e.g., $\overline{M_g}$).
In particular, with the aid of Simpson's approaches to the moduli of coherent
sheaves
over projective schemes (\cite{Simpson94}), our  construction of the universal
moduli space $$M_g({\cal O}, P) \longrightarrow \overline{M_g}$$
is much shorter than the approaches for the similar
moduli problems in \cite{Caporaso94} and \cite{Pandharipande94a}.
See \S\S 8 and 9.
\end{say}

\begin{say} {\sl Koll\'ar's Approaches (as opposed to the traditional  GIT
\cite{GIT})}.
The solution of Question RGIT resembles Koll\'ar's approaches to algebraic
quotient spaces
(cf. Conjecture 1.1 of  \cite{Kollar95a}).
Now, one should ask how to do all of the above via Koll\'ar's approaches
(\cite{Kollar95a} and \cite{Kollar90}).
Theorem 2.14 of \cite{Kollar95a} sheds lights
on this question. But the insistence on the \'etaleness of the the equivariant
map
is too restrictive for universal moduli. There are something
more on the relative quotients that one can do with Koll\'ar's approaches and
in many cases
Koll\'ar's approaches are easier to apply than GIT. In light of the fact that
many moduli spaces of higher dimensional varieties have no
satisfactory GIT constructions, some alternatives are necessary in order to
build universal moduli spaces effectively (for example, over the moduli space
of
surfaces of general type \cite{Kollar94}).
This quest even includes the case of $\overline{M_{g,n}}$
which so far has no (satisfactory) GIT construction.
We will return to these topics  in an upcoming paper \cite{Hu95}.
\end{say}

\medskip\noindent
{\bf Acknowledgements} The prototype of RGIT has appeared in \cite{GIT},
\cite{Simpson94},
\cite{Reichstein89},  and \cite{Pandharipande94a}. The author happily
acknowledges
some helpful conversations and/or correspondence
 with Aaron Bertram, Laurance Ein, David Gieseker, J\'anos Koll\'ar,
Eyal Markman on quotients and/or moduli spaces,
and useful comments from Lev Borisov.  Thanks are due to  Robert Lazarsfeld
for informing him of the references on limit linear series.
 He also thanks  John Millson for  many interesting
communications on moduli spaces of polygons. Finally,
he appreciates the encouragement and support from Dan Burns, Ching-Li Chai,
Igor Dolgachev, Victor Guillemin, Jun Li, Ralf Spatzier, Alex Uribe,
and S.T. Yau during the course of this work.

\bigskip

\begin{center}
{\bf I. Relative GIT} % and Moment Maps}
\end{center}

\section{Equivariant theory of moment maps}

\begin{say}
A {\sl moment map} for a symplectic action of a compact
Lie group $K$ on a symplectic manifold $(M, \omega)$ with a  symplectic form
$\omega$
is a smooth map $$\Phi:M \to {\frak k}^*$$
where ${\frak k}^*$ is the linear dual of the Lie
algebra ${\frak k}$ of $K$,  satisfying the following two properties:
\begin{enumerate}
\item $\Phi$ is equivariant with respect to the given action of $K$ on $M$ and
the
co-adjoint action of $K$ on ${\frak k}^*$;
\item for any $a \in {\frak k}$, $d\Phi(\xi) \cdot a = \omega(\xi, \xi_a)$ for
all
vector fields $\xi$ on $M$,  where
$\xi_a$ is the vector field generated by $a$.
%$\imath_{\omega}(\xi^\sharp) = d(\xi\circ \Phi)$.
\end{enumerate}
The equation in (ii) determines
the moment map by an additive constant in the set ${\frak z}^*$
of the {\sl central elements} of ${\frak k}^*$ (i.e.,
${\frak z}^*$ = the set of invariants of the $K$-coadjoint action = the linear
dual
of the Lie algebra ${\frak z}$ of the
center group of $K$).
\end{say}

When a moment map exists,
we call the the symplectic action a Hamiltonian action
and the symplectic form $\omega$  a ($K$-) Hamiltonian symplectic form.

 Fix a maximal torus $T_k$ of $K$. We shall use $T$ to denote
the complexification of $T_k$. Let $\hbar^*_+$ is a positive
Weyl chamber in the linear dual
 $\hbar^*$ of the Lie algebra $\hbar$ of the fixed maximal torus $T_k$ of $K$.
Since $\hbar^*_+$ parametrizes the orbit space of the coadjoint action of $K$
on ${\frak k}^*$,
we obtain,  by quotienting the $K$-adjoint action,
 the so-called reduced moment map $$\Phi_{\text{red}}: M \to \hbar^*_+.$$

\begin{say}
Although a Hamiltonian symplectic form $\omega$  only determines the moment map
up to an additive central element in ${\frak k}^*$.
However, as pointed out by Atiyah \cite{Atiyah82},
there is always a {\sl canonical choice},  $\Phi^\omega$, the one such that
$$\int_M \Phi^\omega \omega^{\frac{1}{2}\dim M}= 0. $$
Such a moment map will be called the canonical moment map associated to
$\omega$.
All other moment maps defined by $\omega$ have the form
$\Phi^{\widetilde{\omega}} = \Phi^\omega - \mu$ where $\widetilde{\omega}
=(\omega, \mu)$
and $\mu \in {\frak z}$ (consult \ref{say:momentcone} in the sequel).
\end{say}

\begin{rem}
\label{rem:momentmapforsingularspaces}
 For any invariant closed subset $Z$ in $M$ (possibly singular), one can simply
define
the moment map for the action of $K$ on $Z$ to be the restriction of the total
moment map,
for which we shall still use the same notation $\Phi$ if no confusion should
emerge.
\end{rem}

\begin{say} Given any value $p$ of a  moment map $\Phi$,
the orbit space $\Phi^{-1}(\text{\sl O}_p)/K$, which is called the symplectic
reduction at $p$,
 carries a natural induced symplectic structure away
from singularities (\cite{MarsdenWeinstein}).
 Here $\text{\sl O}_p$ is the coadjoint orbit through $p$.
\end{say}

\begin{say}
Now we shall consider equivariant maps and
the relations between their induced moment maps.

%Although in general, (truly) rational equivariant morphisms
%should be considered  for relative GIT, to simplify the exposition,
%we stick ourselves with regular morphisms.

Let  $Y$ and $X$ be smooth projective  complex algebraic varieties acted on by
complex reductive
algbraic group $G'$ and $G$ respectively.   Consider an algebraic surjection
$$\pi: Y \longrightarrow X$$
that is equivariant with respect to a group epimorphism (surjective
homomorphism)
$$\rho: G' \longrightarrow G,$$
that is,
$$\pi (g' \cdot y) = \rho (g') \cdot \pi (y) \; \text{for all} \; y \in Y \;
\text{and}\;
g' \in G'. $$

One can choose a maximal subgroup $K$ of $G$ and
a maximal subgroup $K'$ of $G'$ so that $\rho$ restricts to
a group homomorphism (still denoted by $\rho$)
$\rho: K' \longrightarrow K$ and $\rho: G' \longrightarrow G$ is the
complexification
of $\rho: K' \longrightarrow K$.
\end{say}

\begin{say}
Symplectic (K\"ahler) forms and  line bundles on $Y$ will be denoted by letters
$\eta$ and $M$, while symplectic (K\"ahler) forms and  line bundles on $X$ will
be
denoted by $\omega$ and $L$.
\end{say}

\begin{say}
\label{limitofsymplecticforms}
Let $\omega$ be a symplectic (K\"ahler) form on $X$.
Then $\eta_0 = \pi^* \omega$ is a closed two form on $Y$.
We'd like to define a moment map for this degenerated two form $\eta_0$.
The trick is to use nearby moment maps and then take the limit.
 Let
$\eta (t), t \in (0, \epsilon]$ %\rightarrow {\frak K}^G(Y)$
be a continuous path of symplectic (K\"ahler) forms
on $Y$ such that $\eta_0 = \lim_{t \to 0} \eta (t)$.
Here $\epsilon$ is a small positive number.  (This can be done because
$\eta_0 = \pi^* \omega$ lies on the boundary of the K\"ahler cone for $Y$.)
Define  $$ \Phi^{\eta_0} (y) = \lim_{t \to 0}  \Phi^{\eta (t)} (y)$$
where $ \Phi^{\eta (t)}$ are the canonical moment maps for $\eta (t)$.
\end{say}

\begin{prop}
\label{prop:onlimitforms}
Keep the notations of \ref{limitofsymplecticforms}. Then we have
\begin{enumerate}
\item $\Phi^{\eta_0}$ is a $K'$-equivariant differentiable map;
\item $\Phi^{\eta_0}$ satisfies the differential equation
$$ d\Phi^{\eta_0} (\xi) \cdot a = \eta_0 (\xi, \xi^Y_a)$$
for every vector field $\xi \in TY$ and $a \in {\frak k}'$, where $\xi^Y_a$ is
the
vector field on $Y$ generated by the element $a$;
\item $\Phi^{\eta_0}$ is constant on every fiber of $\pi$.
\end{enumerate}
\end{prop}

\begin{pf}
 The equivariancy of the map  $\Phi^{\eta_0}$ is obvious.
Now since $X \times [0, \epsilon]$ is compact, we have that $\Phi^{\eta (t)}$
converges uniformly
as $t \to 0$. Thus $ \Phi^{\eta_0} = \lim_{t \to 0}  \Phi^{\eta (t)}$ is
differentiable.
This proves (i).
Next,  given any fixed  vector field $\xi \in TY$ and fixed $a \in {\frak k}'$,
we have
$$ d\Phi^{\eta (t)} (\xi) \cdot a = \eta (t) (\xi, \xi_a)$$
for all $t \in (0, \epsilon]$. Passing to the limit as $t \to 0$, we obtain
 $$ d\Phi^{\eta_0} (\xi) \cdot a = \eta_0 (\xi, \xi_a).$$ This proves (ii).
To show (iii), picking any $\xi \in T \pi^{-1}(x) \subset TY$ for any $x \in
X$.
Recall that $\eta_0 = \pi^* \omega$. Now we must have
$$d \Phi^{\eta_0} (\xi) \cdot a = \eta_0 (\xi, \xi_a) = \pi^* \omega (\xi,
\xi_a)
=\omega (d\pi (\xi), d\pi (\xi_a)) = 0 $$
since $d\pi (\xi) = 0$. Because $a$ is an arbitrary element in ${\frak k}'$,
we obtain $d \Phi^{\eta_0} (\xi) = 0$. This implies that $\Phi^{\eta_0}$ is
constant
along every fiber of $\pi$. The proposition is thus proved.
\end{pf}

\begin{rem}
The map $\Phi^{\eta_0}$ may be considered as a moment map defined by
the (pre-symplectic) form $\eta_0$. As the limit of some (true) moment maps,
it enjoys a number of properties of the usual momentum mapping such as
the convexity.  In general, however, a moment map defined by
a  pre-symplectic form may have  much more complicated image that needs not to
be
convex (see \cite{KarshonTollman93}).
\end{rem}

\begin{say}
\label{say:killingmetrics}
Choose a $K$-equivariant metric on ${\frak k}$ and a $K'$-equivariant metric on
${\frak k'}$
respectively so that $d\rho$ preserves the metrics. Then these two metrics lead
two natural
isomorphisms ${\frak k} \cong {\frak k}^*$ and  ${\frak k'} \cong {\frak
k'}^*$, making
the following diagram commutes
\newpage
%\begin{equation}
%\begin{CD}
%{\frak k'}^* @>{\cong}>> {\frak k'}  \\
%@VV{(d\rho)^*}V @VV{d\rho}V \\
%{\frak k}^* @>{\cong}>> {\frak k}
%\end{CD}
%\end{equation}
$${\frak k'}^*\;\;\;\;\;\;\;\;{\buildrel {\cong}\over\longrightarrow}
\;\;\;{\frak k'} $$
$$\;\;\; \uparrow \! (d\rho)^* \;\;\;\;\; \;\;\;\downarrow \! d\rho$$
$${\frak k}^* \;\;\;\;\;\;\;\;{\buildrel {\cong}\over\longrightarrow}
\;\;\;{\frak k}$$
where $d\rho$ is the differential of the map $\rho$
and $(d\rho)^*$ is the codifferential of the map $\rho$, i.e., the linear
dual of the differential $d\rho$.
\end{say}

\begin{thm}
\label{thm:descendingmomentmaps}
The map $\Phi^{\eta_0}$ descends to the moment map
$\Phi^{\omega}$ for the action of $K$ on $X$ with respect to the symplectic
form $\omega$.
In particular, we have the following commutative diagram
$$Y \;\;\;{\buildrel {\Phi^{\eta_0}}\over\longrightarrow} \;\;\;{\frak k'}^*
\;\;\;\;\;\;\;\;
{\buildrel {\cong}\over\longrightarrow} \;\;\;{\frak k'} $$
$$\;\;\;\downarrow \!{\pi} \;\; {\buildrel {\Psi^{\eta_0}}\over\nearrow} \;\;
\uparrow \! (d\rho)^*
\;\;\;\;\; \;\;\;\downarrow \! d\rho$$
$$X\;\;\;{\buildrel{\Phi^{\omega}}\over\longrightarrow} \;\;\;{\frak k}^*
\;\;\;\;\;\;\;\;
{\buildrel {\cong}\over\longrightarrow} \;\;\;{\frak k}$$
\end{thm}

\begin{pf}
By Proposition \ref{prop:onlimitforms} (3),
 the map $\Phi^{\eta_0}$ desends to a well-defined differentiable map
$\Psi^{\eta_0}$ from $X$ to ${\frak k'}^*$ such that $\Phi^{\eta_0}
=\Psi^{\eta_0} \circ \pi$.
Define $\Phi^{\omega}$ to be the composition:
$$\Phi^{\omega}:X \;\;\stackrel{\Psi^{\eta_0}}{\longrightarrow}\;\;
{\frak k'}^* \;\; \stackrel{\cong}{\longrightarrow}\;\; {\frak k'}\;\;
\stackrel{d\rho}{\longrightarrow} \;\;{\frak k}\;\;
 \stackrel{\cong}{\longrightarrow}\;\; {\frak k}^*.$$
Clearly, this is a differentiable map, and  we have that $\Psi^{\eta_0} =
(d\rho)^* \Phi^\omega$.

To show that $\Phi^{\omega}$ is equivariant with respect to the given
$K$-action on $X$
and the coadjoint action on ${\frak k}^*$,
take any element $k' \in K'$ and $y \in Y$, we then have
$$\Phi^{\omega} (\rho(k') \cdot \pi(y)) = \Phi^{\omega} (\pi (k' \cdot y))
=  d\rho \circ \Psi^{\eta_0} (\pi (k' \cdot y)) $$
$$= d\rho \Phi^{\eta_0}(k' \cdot y) = d\rho \text{Ad}(k') \cdot
\Phi^{\eta_0}(y)
= \text{Ad}(\rho(k'))  d\rho \Phi^{\eta_0}(y) $$
$$=\text{Ad}(\rho(k'))  d\rho \Psi^{\eta_0} (\pi (y))
= \text{Ad}(\rho(k')) \Phi^{\omega} (\pi (y)).$$
Here we have used the identity $d\rho \text{Ad}(k') = \text{Ad}(\rho(k'))
d\rho$
coming from the equivariancy of the homomorphism $\rho$.
Using the fact that both $\pi$ and $\rho$ are surjective, we obtain
$$\Phi^{\omega} (k \cdot x) = \text{Ad}(k) \Phi^{\omega} (x), \; \text{for all}
\; k \in K \;
\text{and}\; x \in X.$$

To check that it satisfies the differential equation in the definition
of a moment map (1.1 (2)),
notice that $d\pi :TY \rightarrow TX$ is surjective everywhere but on a lower
dimensional
locus and $d\pi (\xi^Y_{a'}) = \xi^X_{d\rho (a')}$
for any $a' \in {\frak k'}$ because $\pi$ is $\rho$-equivariant, where again,
$\xi^Y_{a'}$ is the vector field on $Y$ generated by $a' \in {\frak k'}$,
while $\xi^X_{d\rho (a')}$ is the vector field on $X$ generated by $d\rho (a')
\in {\frak k}$.
Now for any $\xi \in TY$ and $a' \in {\frak k'}$, we have
$$\omega (d\pi \xi, \xi^X_{d\rho (a')}) = \omega (d\pi \xi, d\pi (\xi^Y_{a'}))
=\pi^* \omega (\xi, \xi^Y_{a'})  = \eta_0(\xi, \xi^Y_{a'}) $$
$$= d\Phi^{\eta_0} (\xi) \cdot a' \; \;\;(\text{because of
\ref{prop:onlimitforms} (2)})$$
$$=d \Psi^{\eta_0} d \pi (\xi) \cdot a' \;\;\; (\text{because}\;
\Phi^{\eta_0} = \Psi^{\eta_0} \circ \pi)$$
$$= (d \rho)^*  d\Phi^{\omega} d \pi (\xi) \cdot  a' \;\;\;
(\text{because}\; \Psi^{\eta_0} = (d \rho)^* \Phi^{\omega})$$
$$=d\Phi^{\omega}(d\pi \xi) \cdot  d \rho (a') \;\;\; (\text{because of
\ref{say:killingmetrics}}).$$
By the surjectivity of $d\rho$, this means that the differential equation
$$d\Phi^{\omega}(\xi^X) \cdot a = \omega (\xi^X, \xi^X_a)$$ holds for
vector fields  $\xi^X \in TX$ almost everywhere.
Thus by continuity, it holds for all $\xi^X$ in $TX$ everywhere.
That is, $\Phi^{\omega}$ is a moment map
for the action of $K$ on $X$ with respect to the symplectic form $\omega$.
Since $\Phi^{\eta (t)}$ are the canonical moment maps for $\eta (t)$, it is
easy to check that
$\Phi^{\omega}$, as the limit of  $\Phi^{\eta (t)}$, is the canonical moment
map for $\omega$.
\end{pf}

%\begin{rem} Suppose now that $\pi: Y \rightarrow X$ is a projective morphism
%%between two
%(possibly singular) quasi-projective varieties and is equivariant with respect
%%to
%the epimorphism $G' \rightarrow G$.  Choose any equivariant projective
%%embeddings
%$$Y \hookrightarrow X \times {\Bbb P}^M  \hookrightarrow {\Bbb P}^N  \times
%%{\Bbb P}^M$$
%$$\downarrow \!{\pi} \;\;\;\;\;\;\;\;\;\; \;\;\;\;\;\;\; \;\;\;\;\;\;\;
%%\downarrow \!\hbox{proj}$$
%$$X \;\;\;\;\;\;\;\;\;\; \;\; \hookrightarrow  \;\;\;\;\;\;\;\;\;\; \;\; {\Bbb
%%P}^N  $$
%Then one can apply Theorem \ref{thm:descendingmomentmaps} to the equivariant
%%projection
%${\Bbb P}^N  \times {\Bbb P}^M \rightarrow {\Bbb P}^N$. By restricting the
%%total moment
%maps to $\pi: Y \rightarrow X$ (cf. Remark \ref{rem:restrictingmomentmaps}),
%%we can extend
% Theorem \ref{thm:descendingmomentmaps} to the equivariant  projective
%%morphism
% $\pi: Y \rightarrow X$ where $X$ and $Y$ may be singular and non-compact. In
%%particular,
%the following commutative diagram still makes sense
%$$Y \;\;\;{\buildrel {\Phi^{\eta_0}}\over\longrightarrow} \;\;\;{\frak k'}^*
%%\;\;\;\;\;\;\;\;
%{\buildrel {\cong}\over\longrightarrow} \;\;\;{\frak k'} $$
%$$\;\;\;\downarrow \!{\pi} \;\; {\buildrel {\Psi^{\eta_0}}\over\nearrow} \;\;
%%\uparrow \! (d\rho)^*
%\;\;\;\;\; \;\;\;\downarrow \! d\rho$$
%$$X\;\;\;{\buildrel{\Phi^{\omega}}\over\longrightarrow} \;\;\;{\frak k}^*
%%\;\;\;\;\;\;\;\;
%{\buildrel {\cong}\over\longrightarrow} \;\;\;{\frak k}$$
%\end{rem}

\begin{say}
\label{deformationofmomentmaps}
Finally, we remark that when $G' = G$,
the above implies there is a deformation of moment maps
$$\Phi^{\eta (t)}: Y \rightarrow {\frak k}^*$$
from $\Phi^{\eta (\epsilon)}: Y \rightarrow  {\frak k}^*$ to
$\Phi^{\omega}: X \rightarrow {\frak k}^*$,
where $\Phi^{\omega} (X) = \Phi^{\eta_0} (Y)$.
When  $\dim G' > \dim G$, the dimension of $\Phi^{\omega} (X)$ is less than
those of  $\Phi^{\eta (t)}(Y) (t \in (0, \epsilon])$. In this case, we say that
$\Phi^{\omega}$ is a degeneration of  $\Phi^{\eta (t)} (t \in (0, \epsilon])$.
\end{say}

\section{$G$-Effective ample cone}
Much of what follows is taken from \cite{DolgachevHu} and \cite{Hu94}.

\begin{say}
In this section we assume that $X$ is a smooth projective complex algebraic
variety\footnote{We point out
that much of results in this section can be extended to
K\"ahler category (i.e., $X$ being
K\"ahler manifolds only).
Since our primary applications will be algebraic moduli spaces, we are content
with working in the category of algebraic varieties.}
 over ${\Bbb C}$ acted on by a reductive group $G$ with a fixed maximal compact
form
$K$.
\end{say}

\begin{thm} \text{(cf. Theorem 2.3.6, \cite{DolgachevHu}.)}
\label{thm:homologicalequivalence} Let $X$ be a  smooth projective variety
acted on by
a reductive complex algebraic group $G$.
Let $\omega$ and $\omega'$ be two $K$-equivariant K\"ahler forms. Suppose that
$\omega$ is cohomological equivalent to  $\omega'$. Then $\Phi^\omega =
\Phi^{\omega'}$.
\end{thm}
\begin{pf} The proof is the same as that for Theorem 2.3.6, \cite{DolgachevHu}.
%By assumption, $\omega' = \omega + \frac{i}{2\pi}d'd'' \hbox{log}(f)$
%for some positive $K$-invariant function $f$. Thus we can write  $\omega' =
%%\omega + d \theta$
%where $\theta$ is a $K$-invariant 1-form of type (0,1).  By definition of the
%%moment map
\end{pf}

\begin{say}  Let $X$ is a smooth projective variety.
We set ${\frak K}^G(X)$ to be the collection of all
$K$-equivariant K\"ahler forms that are compatible with the algebraic
$G$-action modulo cohomological equivalence. Theorem
\ref{thm:homologicalequivalence} says
that there is a well-defined canonical moment map $\Phi^{[\omega]}$ for each
element
$[\omega]$ of ${\frak K}^G(X)$. For notational simplicity,
we shall omit the use the bracket ``[\;\;]''
when it is not likely to cause confusion.
\end{say}

\begin{say}
\label{say:momentcone}
When $X$ is smooth,
set ${\frak M}^G(X)$ to be
the collection of all pairs consisting of an element in ${\frak K}^G(X)$
and a moment map defined by it. Thus
$${\frak M}^G(X) \cong {\frak K}^G(X) \times {\frak z}^*.$$
We shall adopt the following conventional scheme:
{\sl  $\widetilde{\omega}$ denotes an element in
${\frak M}^G(X)$ with its underlying (Hamiltonian)
K\"ahler form symbolized by $\omega \in {\frak K}^G(X)$. Such a symbol
$\widetilde{\omega}$
will be referred as an enriched symplectic K\"ahler form}.
Thus an enriched  symplectic (K\"ahler) form  $\widetilde{\omega}$ in  ${\frak
M}^G(X)$ has the form
$(\omega, \Phi^\omega - \mu)$ or simply $(\omega, \mu)$, where $\mu \in  {\frak
z}^*$.
We frequently write $\Phi^\omega - \mu$ by $\Phi^{\widetilde{\omega}}$.
\end{say}

\begin{rem} All of the results in \S 1 are valid without
modification if the symplectic (K\"ahler) forms ($\omega$ and $\eta (t)$, etc.)
are replaced by the enriched symplectic (K\"ahler) forms
($\widetilde{\omega}$ and $\widetilde{\eta} (t)$, etc.)
\end{rem}

\begin{rem}
If $X$ is (possibly) singular,  we  can set  ${\frak M}^G(X)$ to be the
cone spanned by the images of all linearized ample line bundles in
$\hbox{NS}^G(X)\otimes_{\Bbb Z}
{\Bbb R}$ (\cite{DolgachevHu}).
%Nevertheless, in spite of a little abusing terms, we shall  call a general
%%element
%in  ${\frak K}^G(X)$ a K\"ahler form (after all it is a  K\"ahler form away
%%from singularities
%of $X$ up to cohomological equivalence) and an integral element a Hodge form.
\end{rem}

\begin{defn} Assume that  $X$ is smooth. The $G$-effective ample  cone
is a subcone of ${\frak M}^G(X)$ defined as follows:
$${\frak E}^G(X) =
\{(\omega, \mu) \in {\frak  M}^G(X)\; |\; \mu \in \Phi^\omega(X) \};$$
\end{defn}

\begin{rem} When $X$ is  (possibly) singular, we define
${\frak E}^G(X)$ to be the subcone of  ${\frak M}^G(X)$ that is spanned by
the images of $G$-effective ample line bundles (that is, spanned by the ones
such that
$X^{ss}(L) \ne \emptyset$, see \cite{DolgachevHu}). When $X$ is actually
smooth,
the two definitions are equivalent (\cite{Hu94}).
\end{rem}

\begin{say}
${\frak E}^G(X)$ projects to
${\frak K}^G(X)$ whose fiber at $\omega \in {\frak K}^G(X)$
is the intersection of the moment map image $\Phi^\omega(X)$ with ${\frak
z}^*$.
This fiber  $\Phi^\omega(X) \cap {\frak z}^* = \Phi_{red}^\omega(X) \cap {\frak
z}^*$
is a convex compact polytope.
In fact, it is not hard to see that $\Phi^\omega(X) \cap {\frak z}^*$
is the image of $X$ under the {\it canonical} moment map attached to the
induced action
of the center group of $K$.
\end{say}

\begin{rem}
${\frak E}^G(X)$ may be an empty subset of
${\frak M}^G(X)$ when $G$ is semisimple.  For example, when $X=G/P$ is a
generalized flag
variety. However, ${\frak E}^G(X \times G/B)$ is never empty.
In particular,  ${\frak E}^G(X)$ is never empty when $G$ is a torus.
\end{rem}

\begin{say}
 Recall from \cite{Hu94} that
for any point $x \in X$, we can define a (generalized Hilbert-Mumford)
numerical function
$$M^\bullet (x) : {\frak  E}^G(X) \rightarrow {\Bbb R}$$ whose value
$M^{\widetilde{\omega}}(x)$ at  $ \widetilde{\omega} \in  {\frak  E}^G(X)$
is defined as the {\sl signed} distance
from the origin to the boundary of $\Phi^{\widetilde{\omega}} (\overline{G
\cdot x})$:
it takes a positive value if $0$ is outside of $\Phi^{\widetilde{\omega}}
(\overline{G \cdot x})$;
 it takes a nonpositive value otherwise. Using this numerical function
we have the following criteria for (K\"ahler) semistabilities,
\begin{enumerate}
\item $X^{ss}(\widetilde{\omega})$ = the set of the
 semistable points with respect to $\widetilde{\omega}$ \\
= $\{x \in X| M^{\widetilde{\omega}}(x) \le 0\}$;
\item $X^s(\widetilde{\omega})$ = the set of the
stable points with respect to $\widetilde{\omega}$ \\
= $\{x \in X| M^{\widetilde{\omega}}(x) < 0\}$;
\item $X^{us}(\widetilde{\omega})$ = the set of the unstable (or
non-semistable)
 points with respect to $\widetilde{\omega}$ \\
= $\{x \in X| M^{\widetilde{\omega}}(x) > 0\}$;
\item{(iv)} $X^{sss}(\widetilde{\omega})$ = the set of the
strictly semistable points with respect to $\widetilde{\omega}$ \\
= $X^{ss}(\widetilde{\omega}) \setminus X^s(\widetilde{\omega})$.
\end{enumerate}
\end{say}

\begin{rem}
When a group is needed to be specified, we will write
$X^{ss}_G(\widetilde{\omega})$,
$X^{s}_G(\widetilde{\omega})$, etc. This applies especially when there is a
redutive subgroup
$H$ of $G$. In this case, the restriction map ${\frak M}^G(X) \to {\frak
M}^H(X)$ (the $H$-moment map
is obtained from the $G$-moment map
by the orthogonal projection ${\frak k}^* \to {\frak h}^*$ where ${\frak h}$ is
the
Lie algebra of a suitable compact form of $H$) induces a linearization in
${\frak M}^H(X)$
for each  linearization $\widetilde{\omega}$ in  ${\frak M}^G(X)$. To specify
this effect,
we will write $X^{ss}_H(\widetilde{\omega})$ ($X^{s}_H(\widetilde{\omega})$,
etc) for the set
of semistable (stable, etc) points for the action of $H$.
\end{rem}

By the works of Kempf-Ness (for the algebro-geometric cases \cite{KempfNess78})
and Kirwan
(for the K\"ahler generalizations \cite{Kirwan84}),
there is a Hausdorff quotient topology
on  $X^{ss}(\widetilde{\omega})/\!/G$ such that it  contains
the orbit space $X^s(\widetilde{\omega})/G$ as a dense open subset. It is a
K\"ahler space
and has a K\"ahler form induced from $\omega$ away from the singularities.
Moreover,
it is homeomorphic to the symplectic reduction
$(\Phi^{\widetilde{\omega}})^{-1}(0)/K$.

\begin{thm} \text{(Kempf-Ness-Kirwan)}
\label{thm:Kempf-Ness-Kirwan}
Let $\widetilde{\omega}=(\omega, \mu)$ be an enriched   K\"ahler form in
${\frak E}^G(X)$.
Then  $(\Phi^{\widetilde{\omega}})^{-1}(0) \subset X^{ss}(\widetilde{\omega})$
and the inclusion induces
a homeomorphism $$(\Phi^{\widetilde{\omega}})^{-1}(0)/K {\buildrel
{\cong}\over\longrightarrow}
X^{ss}(\widetilde{\omega})/\!/G.$$
In case that $\widetilde{\omega}$ is integral (i.e., coming from a linearized
ample line bundle),
then $X^{ss}(\widetilde{\omega})/\!/G$ carries a projective structure.
\end{thm}

%\begin{rem} Since ${\frak E}^G(X)$ is convex, hence connected. Theorem
%%\ref{thm:Kempf-Ness}
%implies that

\begin{rem}
The $G$-effective ample cone takes care of symplectic reductions at central
values of ${\frak k}^*$ and identifies them with K\"ahler quotients of $X$ by
the complex reductive group $G$.
To include symplectic reductions at non-central values, one has to consider the
so-called
enlarged moment cone (see \cite{Hu94}) and use the so-called shifting trick to
identify
them with the K\"ahler quotients on $X \times G/P$ by the diagonal action of
$G$.
The relative GIT for the morphisms $X \times G/B \to X \times G/P$ has been
studied
and linked to degenerated quotients of $X \times G/B$ in \cite{Hu94}. So, in
this paper,
we stick with just the $G$-effective ample cone.
\end{rem}

\begin{say}
 The union of the zero sets of $M^\bullet(x): {\frak E}^G(X) \rightarrow {\Bbb
R}$
for all $x$ with isotropy subgroups of positive dimensions equals the union
${\cal W}$ of all
walls in ${\frak E}^G(X)$ (see \cite{DolgachevHu} and \cite{Hu94}).
A connected component of ${\frak E}^G(X) \setminus {\cal W}$
is a chamber. Linearizations in the same chamber define the same notion of
stabilities.
\end{say}

We shall need the following results  in the sequel.

\begin{say}
Given any element $\widetilde{\omega}
 \in {\frak E}^G(X)$, hence a {\it unique} moment map
$\Phi^{\widetilde{\omega}}= \Phi$.
We have a stratification, the momentum Morse stratification
with respect to $\Phi$,  $X = \cup_{\beta \in {\bf B}} S_\beta$ induced by the
norm square
$|\!|\Phi|\!|^2$ of the moment map (\cite{Kirwan84}).
The strata $S_\beta$ and
their indexes $\beta$ can be described  as follows:
$$S_\beta = \{x\in X | \; \beta \  \hbox {is the unique closest point to 0 of
}\ \Phi_{\red}
 (\overline {G\cdot x})\}.$$
\end{say}

\begin{thm}
\label{thm:finiteness} {\rm (\cite{DolgachevHu}, \cite{Hu94})}
There are only finitely many momentum Morse stratifications.
\end{thm}

\begin{defn} {\rm (\cite{Hu94})}
\label{defn:thinmorsestratification}  For any  momentum Morse stratification of
$X$,
we choose precisely one stratum from it. Then
the intersection of all the chosen strata is called
 a {\it thin} momentum Morse stratum provided that it is not empty.
There are only finitely many such strata. The
 {\it thin} momentum Morse strata form a stratification of $X$.
 \end{defn}

\begin{prop}  {\rm (\cite{Hu94})}
Two points in the same {\it thin} momentum Morse stratum
give rise to the same numerical function $M^\bullet (x): {\frak E}^G(X)
\rightarrow {\Bbb R}$.
\label{cor:tmm=snf}
\end{prop}

\begin{pf} The proof is short. So we repeat it here (see \cite{Hu94}).
Let $x$ and $y$ be two points of a {\it thin} momentum Morse stratum.
Then by definition, for any $\widetilde{\omega} \in  {\frak E}^G(X)$,
$x, y$ belong to the same momentum Morse stratum
$S^{\widetilde{\omega}}_{\beta(\widetilde{\omega})}$
for some index  $\beta(\widetilde{\omega})$.
That is, $M^{\widetilde{\omega}}(x) = M^{\widetilde{\omega}}(y)
 = |\!|\beta(\widetilde{\omega})|\!|$ for all $\widetilde{\omega} \in {\frak
E}^G(X)$.
Hence $$M^\bullet(x) = M^\bullet(y) : {\frak E}^G(X) \rightarrow {\Bbb R}.$$
\end{pf}

\begin{cor}
\label{cor:finitenumericalfunctions}
There are only finitely many numerical functions
$$M^\bullet (x): {\frak E}^G(X) \rightarrow {\Bbb R}.$$
(And the points in the same {\it thin} momentum Morse stratum give rise to
 the same numerical function $M^\bullet (x)$.)
\end{cor}

\begin{rem} To close this section we make the following useful remark.
We view the integral points in ${\frak K}^G(X)$ (${\frak E}^G(X)$) as being
induced from
the Hodge metrics of  (linearized) ample line bundles. We shall think of
{\it these integral points and (linearized) ample line bundles interchangably.}
 By default, {\sl the K\"ahler quotiens and maps associated to these integral
points
will be projective as in the traditional GIT cases.}
\end{rem}

\section{A relative GIT: first cases}

\begin{say} In the rest of the paper, unless specified otherwise, %and in
%%general
 we shall assume that
$\pi:Y \rightarrow X$ is a projective morphism between
two (possibly singular) quasi-projective algebaric varieties that is
equivariant
with respect to an epimorphism $\rho: G' \rightarrow G$
between two reductive complex algebraic groups having  fixed maximal
compact subgroups $K'$ and $K$, respectively.
\end{say}

\begin{say}
\label{say:divideandconquer}
Let $G_0$ be the kernel of $\rho$. Then we have
$$\{1\} \rightarrow G_0 \rightarrow G' \rightarrow G \rightarrow \{1\} .$$
That is, up to a finite central extension, we may think of $G' = G_0 \times G$.
Because GIT problems for a finite group is trivial (i.e., the orbit space
is the natural  solution to the quotient problems for a finite group action),
to simplify the exposition, one may well assume that $G \cong G'/G_0$ is a
subgroup of $G'$
and $G' = G_0 G$.
%However, in keeping with the general theme, we shall not assume so.

Due to the above comments, the relative GIT problem for this case can be
divided into two steps:
\begin{enumerate}
\item The fiberwise GIT problem: $G_0$ acts only on the fibers of $\pi: Y \to
X$;
(One may look at this from a slightly different point of view: $Y/X$ is
projective over the
base scheme $X$. $G_0/X$ is the trivial group scheme over $X$. $G_0/X$ acts on
$Y/X$. This
viewpoint is helpful for some universal moduli problems.)
\item The $G$-equivariancy GIT problem: treat $\pi: Y \to X$ as a
$G$-equivariant map alone.
Here we identify the quotient group $G'/G_0$ with the group $G$ by the
isomorphism induced
by the epimorphism $\rho$ and $G$ acts on $Y$ via this identification.
\end{enumerate}
After those have been done, we can then put the two together and
relate some properly chosen $Y$-quotient by $G'$ to a given $X$-quotients by
$G$.
\end{say}

\begin{prop}
\label{prop:fiberwisegit}
Let $M$ be a $G_0$ linearized relatively ample line bundle
with respect to the morphism $\pi: Y \to X$. Let $Y_x$
denote the fiber
of $\pi$ over a point $x \in X$ and $M_x$ the restriction of
$M$ to $Y_x$. Then we have $Y_x^{ss}(M_x) = Y_x \cap Y_{G_0}^{ss}(M)$.
In particular $Y_{G_0}^{ss}(M)= \cup_{x \in X} Y_x^{ss}(M_x)$.
\end{prop}

\begin{pf} It follows from Proposition 1.19 of \cite{GIT}. See also
\cite{Simpson94}.
\end{pf}

Consequently, since the categorical quotient is universal,
one sees that the morphism $\pi: Y^{ss}_{G_0}(M) \to X$ descends to
a map $$\widetilde{\pi}: Y_{G_0}^{ss}(M)/\!/G_0 \to X$$ with fibers
$Y_x^{ss}(M_x)/\!/G_0$.

\begin{prop}
\label{G0quotientisGinvariant}
If $M$ is a $G'$-linearized relatively ample line bundle
with respect to the morphism $\pi: Y \to X$, then
the $G'$-action on $Y$ descends to a $G'/G_0=G$-action on
$Y_{G_0}^{ss}(M)/\!/G_0$ making
the map $\widetilde{\pi}: Y_{G_0}^{ss}(M)/\!/G_0 \to X$ $G$-equivariant.
\end{prop}

\begin{pf}
First we need to show that $Y_{G_0}^{ss}(M)$ is $G'/G_0=G$-invariant (recall
that
the quotient group $G'/G_0$ is identified with the group $G$ by the isomorphism
induced
by the epimorphism $\rho$). Take a point $y \in Y_{G_0}^{ss}(M)$ and an element
$g \in G$.
Let $s \in \Gamma(Y, M^{\otimes n})^{G_0}$ be an $G_0$-invariant section
such that $Y_s:=\{y' \in Y | s(y') \ne 0 \}$ is affine and contains $y$
(see Defintion 1.7 in \S 4 of \cite{GIT}).
Set $$^gs (y') = g s(g^{-1} y')$$ for all $y' \in Y$. We claim that
$^gs \in \Gamma(Y, M^{\otimes n})^{G_0}$. To see this, pick any element $g_0
\in G_0$.
We have $^{g_0} (^gs)  (y') = g_0 (^gs) (g_0^{-1} y') = g_0 g s(g^{-1} g_0^{-1}
y')
= g \hat{g}_0 s (\hat{g}_0^{-1} g^{-1} y') = g   (^{\hat{g}_0})s (g^{-1} y')
= g s (g^{-1} y') = (^gs) (y')$ where $\hat{g}_0$ is an element in $G_0$ such
that
$g_0 g = g \hat{g}_0$ (note that $G_0$ is a normal subgroup). This implies that
 $^{g_0} (^gs) = ^gs$, i.e., $^gs \in \Gamma(Y, M^{\otimes n})^{G_0}$.
Clearly $Y_{^gs} = g Y_s$ is affine and contains $g y$.
That is,  $g y \in Y_{G_0}^{ss}(M)$.

Since $Y_{G_0}^{ss}(M) \subset Y$ factors to $X$ in
a $G_0$-equivariant way ($G_0$ acts trivially on $X$), by the universality of
categorical quotient,
we see that there is a naturally induced morphism $\widetilde{\pi}:
Y_{G_0}^{ss}(M)/\!/G_0 \to X$.
Next it is straightforward to verify that the $G$-action on $Y_{G_0}^{ss}(M)$
passes naturally to the quotient and the $G$-equivariancy of $\pi:
Y_{G_0}^{ss}(M) \to X$
implies the $G$-equivariancy of   $\widetilde{\pi}: Y_{G_0}^{ss}(M)/\!/G_0 \to
X$.
\end{pf}

\begin{say} Now let us consider the $G$-equivariancy GIT problems. That is,
we are going to treat  morphism $\pi: Y \to X$ as a $G$-equivariant map only.
Our finest results in this section lie in the case when $X^{ss}(L) = X^s(L)$.
This is what we shall sometimes refer as ``good cases'' .
%As some results are also easy without this assumption, we from time to time
%may also steer away, for a moment, to the general case  when $X^{ss}(L) \ne
%%X^s(L)$.
\end{say}

\begin{say}
\label{say:setupphieta} We begin with assuming that $\pi: Y \to X$ is a
projective morphism
between two smooth projective varieties.
Let $\widetilde{\omega} = (\omega, \mu)$ be an enriched  $K$-equivariant
K\"ahler form on $X$
and $\widetilde{\eta}_0 = \pi^* \widetilde{\omega} = (\pi^* \omega, \mu)$
be an enriched closed $K$-equivariant two form on $Y$. Let
$\widetilde{\eta}: (0, \epsilon] \rightarrow {\frak E}^G(Y)$ be a continuous
path
%in a chamber (need not be a top one) of
in ${\frak E}^G(Y)$
such that $ \widetilde{\eta}_0 = \lim_{t \to 0} \widetilde{\eta} (t)$.
As in \ref{limitofsymplecticforms},
define $$ \Phi^{\widetilde{\eta}_0} (y) = \lim_{t \to 0}
\Phi^{\widetilde{\eta} (t)} (y).$$
Then we have
\end{say}

\begin{thm}
\label{thm:comparingstabilities}
Let  $Y \rightarrow X$ be  a $G$-equivariant
projective morphism between two smooth projective varieties. Then
there exists $\delta > 0$ such that
\begin{enumerate}
\item If $x = \pi (y)$ is stable in $X$ w.r.t  $\widetilde{\omega}$,
 then $y$ is stable in $Y$ w.r.t. $\widetilde{\eta} (t), t \in (0, \delta]$;
\item If $x = \pi (y)$ is non-semistable in $X$  w.r.t  $\widetilde{\omega}$
then $y$ is non-semistable in $Y$ w.r.t. $\widetilde{\eta} (t), t \in (0,
\delta]$;
\item $Y^{s}(\widetilde{\eta} (t)) \supset \pi^{-1}(X^{s}(\widetilde{\omega}))$
and $Y^{ss}(\widetilde{\eta} (t)) \subset \pi^{-1}(X^{ss}(\widetilde{\omega}))$
for $t \in (0, \delta]$.
\end{enumerate}
\end{thm}

\begin{pf}   Let $Y= \bigcup_{i=1, \cdots, d} S_i$ be the thin momentum Morse
stratification
of $Y$ (see Definition \ref{defn:thinmorsestratification}).
By Corollary \ref{cor:finitenumericalfunctions}, it suffices to consider
some fixed representing elements in $S_i, i=1, \cdots, d$.
For any $1 \le i \le d$, if $\pi (S_i) \subset X^{sss}(\widetilde{\omega})$,
set $\delta_i =1$.

Otherwise,  pick up some representing elements $y_i \in S_i,  1 \le i \le d$.

If $x_i = \pi (y_i)$ is stable in $X$, then
$0 \in \hbox{int}(\Phi^{\widetilde{\omega}}(\overline{G \cdot x_i}))$.
By the remarks in \ref{deformationofmomentmaps},
a small deformation of $\Phi^{\widetilde{\omega}}(\overline{G \cdot x_i})$
should still contains the origin in its interior. That is, there exists
$\delta_i > 0$
such that
$0 \in \hbox{int}(\Phi^{\widetilde{\eta} (t)}(\overline{G \cdot y}))$ for
$t \le \delta_i$. Thus $y_i$ is stable in $Y$  w.r.t.  $\widetilde{\eta} (t)$
for
 $t \le \delta_i$. %But the stability only depends on the chamber, by the
%%choice of
%the path $\widetilde{\eta} (t)$,

If  $x_i = \pi (y_i)$ is non-semistable in $X$, then
 $0$ is outside of $\Phi^{\widetilde{\omega}}(\overline{G \cdot x_i})$.
By the similar deformation argument as above,
the same is true for its small deformations.
This is,, there exists $\delta_i > 0$
such that
$0 \notin \hbox{int}(\Phi^{\widetilde{\eta} (t)}(\overline{G \cdot y}))$ for
$t \le \delta_i$. In other words,
$y$ is non-semistable in $Y$ w.r.t $\widetilde{\eta} (t)$ for
 $t \le \delta_i$.

Now choose $\delta= \hbox{min}\{\delta_1, \cdots, \delta_d, \epsilon \}$,
the above implies both (1) and (2).

(3) follows from (2).
\end{pf}

\begin{rem}
When $x=\pi (y)$ is strictly semistable, i.e., the orgin is contained
on the boundary of $\Phi^{\widetilde{\omega}}(\overline{G \cdot x})$,
the stability of $y$ depends on the direction of deformation. So in general,
all
 three possible situations
(stable, strictly semistable, non-semistable) may happen.
\end{rem}

\begin{say}
\label{say:traditional GIT}
In keeping with the theme of the traditional GIT,
consider the case that $\widetilde{\omega}$
 is an integral form induced from an ample linearized line bundle
$L$ on $X$. The pullback $\pi^* L$ is only a {\sl nef}
linearized line bundle on $Y$ inducing the 2-form
$\widetilde{\eta}_0 = \pi^* \widetilde{\omega}$. To get an ample linearized
line bundle on $Y$,
we need to choose an arbitrary
{\it relatively} ample linearized  line bundle $M$ on $Y$ and take
 a sufficiently large tensor power of $\pi^* L$. That is,
$\pi^* L^n \otimes M$ is ample for $n \gg 0$ (\cite{EGA}).
Set fractional linearizations $M_n = \frac{1}{n} (\pi^* L^n \otimes M)$.
Clearly, $\lim_{n \to \infty} M_n = \pi^* L$. The purpose of
this scheme is two fold: to get (fractional multiples of) {\it ample}
 linearized line bundles $M_n$ on $Y$;
and to make $\pi^* L$ and $M_n = \frac{1}{n}(\pi^* L^n \otimes M)$
sufficiently close for sufficiently large $n$.
% In other words, for good choice of $L$
%(e.g., $X^{sss}(L) = \emptyset$), $M_n$ ($n \gg 0$) will lie in a top chamber
%%of ${\frak E}^G(Y)$
%that is independent of $M$ and contains  $\pi^* L^n$ (for all large $n$) in
%%its closure.
\end{say}

We need the following easy lemma
\begin{lem}
\label{lem:restriction}
 Let $G$ act on $X$ whose acton is linearized by $L$ and
$i: Z \hookrightarrow X$ be a $G$-invariant closed embedding linearized by
$i^* L$. Then $Z^{ss}(i^* L) = Z \cap X^{ss}(L)$ and $Z^{s}(i^* L) = Z \cap
X^{s}(L)$.
\end{lem}
\begin{pf} This follows from Proposition 1.19 of \cite{GIT}.
\end{pf}

\begin{thm} %\text{(Generalized Reichstein's Theorem)}
\label{thm:generalizedReichstein}
Let $\pi: Y \rightarrow X$ be a $G$-quivariant projective morphism between two
(possibly singular) quasi-projective varieties. Given any linearized ample line
bundle
on $L$ on $X$, choose a {\it relatively} ample linearized  line bundle $M$ on
$Y$. Then there
exists $n_0$ such that when $n \ge n_0$, we have
\begin{enumerate}
\item  $Y^{ss}(\pi^*L^n \otimes M) \subset \pi^{-1}(X^{ss}(L))$;
%If $x = \pi (y)$ is stable in $X$ then $y$ is stable in $Y$;
\item   $Y^{s}(\pi^*L^n \otimes M) \supset \pi^{-1}(X^{s}(L))$
 %If $x = \pi (y)$ is non-semistable in $X$ then $y$ is non-semistable in $Y$;

If in addition, $X^{ss}(L) = X^s(L)$, then
\item $Y^{ss}(\pi^*L^n \otimes M)=Y^{s}(\pi^*L^n \otimes M)=\pi^{-1}(X^{s}(L))=
\pi^{-1}(X^{ss}(L))$. In particular, $\pi^*L^n \otimes M$ lie in the same
chamber
of ${\frak E}^G(Y)$ for all  $n \ge n_0$.
\end{enumerate}
\end{thm}
\begin{pf}
Given a linearized  ample line bundle $L$ on $X$, and
a relatively linearized ample line bundle $M$ for $\pi: Y \rightarrow X$,
by \cite{EGA},
there exists $m_0$ such that when $m \ge m_0$, $L^m$ and
$\pi^* L^{\otimes m} \otimes M$ are very ample.
Consider some  equivariant projective embeddings induced from
$L^m$ and $\pi^* L^{\otimes m} \otimes M$ (for some $m \ge m_0$)
$$Y \hookrightarrow X \times {\Bbb P}^{N'}  \hookrightarrow {\Bbb P}^N  \times
{\Bbb P}^{N'}$$
$$\;\;\;\downarrow \!{\pi} \;\;\;\;\;\;\;\;\;\;\;\;\;
\;\;\;\;\;\;\; \;\;\;\;\;\;\; \downarrow \!\hbox{proj}$$
$$X \;\;\;\;\;\;\;\;\;\; \;\; \hookrightarrow  \;\;\;\;\;\;\;\;\;\; \;\; {\Bbb
P}^N  $$
We may think $L^m$ as the pull back of a very ample line bundle ${\cal
O}_{{\Bbb P}^N}(1)$.
By \ref{say:traditional GIT}, there exists $d_0 > 0$ such that
when $d \ge d_0$, Theorem \ref{thm:comparingstabilities}
 can be applied to the $G$-equivariant projection
${\Bbb P}^N  \times {\Bbb P}^{N'} \rightarrow {\Bbb P}^N$ with respect to the
linearizations ${\cal O}_{{\Bbb P}^N}(1)$ on ${\Bbb P}^N$ and
$\frac{1}{d}({\cal O}_{{\Bbb P}^N}(d) \otimes {\cal O}_{{\Bbb P}^{N'}}(1))$
on $ {\Bbb P}^N  \times {\Bbb P}^{N'}$, respectively . Now restricting
everything
to $\pi: Y \rightarrow X$ (see Lemma \ref{lem:restriction}),
 %(cf. Remark \ref{rem:restrictingmomentmaps}),
 we see that
$\pi: Y \rightarrow X$ %with respect to the linearizations
shares all the properties as in Theorem \ref{thm:comparingstabilities}. This
proves (1) and (2).

All of (3) follows immediately from (1) and (2).
\end{pf}

\begin{rem}
In the case that $M$ is ample and   $X$ and $Y$ are projective,
Theorem \ref{thm:generalizedReichstein}
is largely due to Z. Reichstein (\cite{Reichstein89}) whose proof is completely
different.
But we only need to assume that $M$ is {\it relatively} ample and   $X$ and $Y$
are
quasi-projective. We emphasize the practical convenience that can
result from the {\it relative} ampleness of the line bundle $M$,
because in constructing
moduli spaces,  relative projective embeddings of  relative
Hilbert schemes have been constructed very explicitly by Grothendieck
(hence are ready to use),  whereas  absolute projective  embeddings will not be
offhand.
\end{rem}

\begin{thm}
\label{thm:inducedmapfromG-equivariancy}
Keep the assumption as in Theorem \ref{thm:generalizedReichstein}.
Then the inclusion $Y^{ss}(\pi^* L^n \otimes M) \subset
\pi^{-1}(X^{ss}(L))$ induces a projective morphism
$$\hat{\pi}: Y^{ss}(\pi^* L^n \otimes M)/\!/G \longrightarrow
X^{ss}(L)/\!/G$$ such that
\begin{enumerate}
\item $\hat{\pi}^{-1}([G \cdot x]) \cong \pi^{-1}(x)/G_x$ for any $x \in
X^s(L)$;
\item if $\pi$ is a fibration, $X^{ss}(L) =X^s(L)$,  and
$G$ acts freely on $X^s(L)$, then $\hat{\pi}$ is also a fibration with the same
fibers.
\end{enumerate}
\end{thm}

\begin{pf} The morphism
$\pi$ restricts to give a morphism $Y^{ss}(\pi^* L^n \otimes M) \rightarrow
X^{ss}(L)$.
This in turn induces an obvious morphism
$$Y^{ss}(\pi^* L^n \otimes M)\longrightarrow X^{ss}(L)
 \longrightarrow X^{ss}(L)/\!/G.$$
Now by the universality of the categorical quotient, we get our desired induced
morphism
$$\hat{\pi}: Y^{ss}(\pi^* L^n \otimes M)/\!/G \longrightarrow
X^{ss}(L)/\!/G.$$

To show (1), let $[G \cdot y_1]$ and $[G \cdot y_2]$ be arbitrary two points in
$Y^{ss}(\pi^* L^n \otimes M)/\!/G$ that are mapped down to $[G \cdot x] \in
X^s (L)/\!/G$. Then we must have that $\pi (y_1), \pi (y_2) \in X^s (L)$.
By Theorem \ref{thm:generalizedReichstein} (2), $y_1, y_2 \in Y^{s}(\pi^* L^n
\otimes M)$.
Thus $\pi (G \cdot y_1) = \pi (G \cdot y_2) = G \cdot x$. Hence by applying
some elements of $G$, we may assume that
$y_1, y_2 \in \pi^{-1}(x)$. But, in this case,  $[G \cdot y_1]=[G \cdot y_2]$
if and only if $y_2 = g\cdot y_1$ for some $g \in G$. This implies that $x =g
\cdot x$.
That is, $g \in G_x$. Hence we obtain that
 $\hat{\pi}^{-1}([G \cdot x]) \cong \pi^{-1}(x)/G_x$.

(2) follows from (1) immediately.
\end{pf}

\section{A relative GIT: general cases}

Now we begin to investigate the relation between $G'$-stability on $Y$ and
the $G$-stability on $X$.

\begin{defn} A linearized pair $(L, M)$ for the morphism $\pi: Y \to X$
\label{defn:linearizedpair}
 consists of a $G$-linearized ample line bundle $L$ over $X$ and a
$G'$-linearized
$\pi$-ample line bundle $M$ over $Y$.
\end{defn}

Out of a linearized pair $(L, M)$, we can have two Zariski open subsets
$Y^{ss}_{G_0}(M)$ and $\pi^{-1}(X^{ss}(L))$.  Proposition
\ref{G0quotientisGinvariant}
says that $Y^{ss}_{G_0}(M)$ is $G$-invariant and hence $G'$-invariant.
Obviously  $\pi^{-1}(X^{ss}(L))$ is $G_0$-invariant, hence also $G'$-invariant.

\begin{lem}
\label{triviallinearization}
Let  $(L, M)$ be a linearized pair for the morphism $\pi: Y \to X$.
Then
\begin{enumerate}
\item $Y^{ss}_{G_0}(\pi^* L^n \otimes M) = Y^{ss}_{G_0}(M)$;
\item  $Y^{s}_{G_0}(\pi^* L^n \otimes M) = Y^{s}_{G_0}(M)$;
\end{enumerate}
\end{lem}
\begin{pf}
Notice that $G_0$ preserves the fibers of the morphism $\pi$. Thus we are in
the position
to apply Proposition \ref{prop:fiberwisegit}. Given any point $x \in X$.
Since the restricted line bundle $\pi^* L^n|_{Y_x}$ is trivial, we obtain
$(Y_x)_{G_0}^{ss}((\pi^* L^n \otimes M)|_{Y_x}) = (Y_x)^{ss}_{G_0}(M)$.
Now Proposition \ref{prop:fiberwisegit} implies that
$Y^{ss}_{G_0}(\pi^* L^n \otimes M) = Y^{ss}_{G_0}(M)$.

(2) follows from (1) because of the fact: a point is stable if and only if its
isotropy is
finite and its orbit is closed in the semistable locus.
\end{pf}

Again by 4.6.13 (ii) of \cite{EGA},
there exists an positive integer $n_0$ such that
$\pi^* L^n \otimes M$ is ample when $n \ge n_0$.

\begin{lem}
\label{G'=GG_0}
Assume that $n \ge n_0$. Then
\begin{enumerate}
\item $Y^{ss}(\pi^* L^n \otimes M) = Y^{ss}_{G_0}(\pi^* L^n \otimes M) \cap
 Y^{ss}_{G}(\pi^* L^n \otimes M).$
\item $Y^{s}(\pi^* L^n \otimes M) = Y^{s}_{G_0}(\pi^* L^n \otimes M) \cap
 Y^{s}_{G}(\pi^* L^n \otimes M).$
\end{enumerate}
\end{lem}
\begin{pf} By general nonsense, we have
$$Y^{ss}(\pi^* L^n \otimes M) \subset Y^{ss}_{G_0}(\pi^* L^n \otimes M) \cap
 Y^{ss}_{G}(\pi^* L^n \otimes M).$$
To show the other way inclusion, notice that
$G'=G_0 G$ and any 1-PS of $G'$ can be written as $\lambda \lambda_0$ where
$\lambda$ is a 1-PS of $G$ and $\lambda_0$ is is a 1-PS of $G_0$.
Take a point $y \in Y^{ss}_{G_0}(\pi^* L^n \otimes M) \cap
 Y^{ss}_{G}(\pi^* L^n \otimes M)$. Now by Corollary 2.15 in \S 3 of \cite{GIT},
$$\mu^{\pi^* L^n \otimes M}(y, \lambda \lambda_0)
\le \mu^{M \otimes \pi^* L^{n}}(y, \lambda) + \mu^{\pi^* L^n \otimes M}(y,
\lambda_0) \le 0+0=0.$$
That is, $Y^{ss}(\pi^* L^n \otimes M) \supset Y^{ss}_{G_0}(\pi^* L^n \otimes M)
\cap
 Y^{ss}_{G}(\pi^* L^n \otimes M)$. This shows (1).

To show (2), again by general nonsense,
$$Y^{s}(\pi^* L^n \otimes M) \subset Y^{s}_{G_0}(\pi^* L^n \otimes M) \cap
 Y^s_{G}(\pi^* L^n \otimes M).$$
Similar to the semistable case,
the other way inclusion follows from the fact that $G'=GG_0$
and Corollary 2.15 in \S 3 of \cite{GIT}.
\end{pf}

Now we are in the position to state our main theorem on the comparisons of
$G'$-semistabilities on $Y$ and $G$-semistabilities on $X$.

\begin{thm}
\label{thm:generalcomparisonofstabilities}
Let $(L, M)$ be a linearized pair for the morphism $\pi: Y \to X$.
%Assume that $X^{ss}(L)=X^s(L)$.
Then for sufficiently large $n$
\begin{enumerate}
\item  $Y^{ss}(\pi^* L^n \otimes M) \subset Y^{ss}_{G_0}(M) \cap
\pi^{-1}(X^{ss}(L))$.
% A point from this set is called semistable with respect to   $(L, M)$.
\item  $Y^s(\pi^* L^n \otimes M) \supset Y^s_{G_0}(M) \cap \pi^{-1}(X^s(L))$ \\
%\cup  (Y^{s}_{G_0}(M) \cap \pi^{-1}(X^{ss}(L)))$.
%A point from this set is called semistable with respect to   $(L, M)$.
%\item $Y^{sss}(L, M) = Y^{ss}(L, M) \setminus Y^s(L, M)$. A point from this
%%set is
%called strictly semistable with respect to   $(L, M)$.
%\item  $Y^{us}(L, M)= Y \setminus Y^{ss}(L, M)$.  A point from this set is
%called unstable (or more correctly, non-semistable).
Assume in addition that $X^{ss}(L)=X^s(L)$, then
\item $Y^{ss}(\pi^* L^n \otimes M) = Y^{ss}_{G_0}(M) \cap \pi^{-1}(X^s(L))$.
\item $Y^s(\pi^* L^n \otimes M) =   Y^{s}_{G_0}(M) \cap \pi^{-1}(X^s(L))$.
\end{enumerate}
\end{thm}
\begin{pf}
{}From Lemma \ref{G'=GG_0} (1),
 $$Y^{ss}(\pi^* L^n \otimes M) = Y^{ss}_{G_0}(\pi^* L^n \otimes M) \cap
 Y^{ss}_{G}(\pi^* L^n \otimes M).$$
Then by Lemma \ref{triviallinearization} (1) and
Theorem \ref{thm:generalizedReichstein} (1)
(see also Theorem \ref{thm:comparingstabilities} (3)),
we have
$$Y^{ss}(\pi^* L^n \otimes M) \subset Y_{G_0}^{ss}(M) \cap
\pi^{-1}(X^{ss}(L)).$$
This shows (1), as desired.

(2) also follows in a similar way from
the combination of  Lemma \ref{G'=GG_0} (2), Lemma \ref{triviallinearization}
(2),  and
Theorem \ref{thm:generalizedReichstein} (2).

Now assume that $X^{ss}(L) =X^s (L)$. Then, similar to the above, (3) and (4)
follow directly from the combination of  Lemma \ref{triviallinearization},
Lemma \ref{G'=GG_0}, and Theorem \ref{thm:generalizedReichstein}.
\end{pf}

\begin{rem}
When  $\pi: Y \to \text{\{point\}}$ is the total contraction and $\rho$ is the
homomorphism
$G \to \{1\}$, we have $L = \text{\sl O}_{\text{\{point\}}}$. In this case, the
relative
$\pi$-ampleness is equivalent to the absolute ampleness and the effect of
the  a  linearized pair $(L, M)$ is equivalent to that of the  linearized ample
line bundle $M$. Thus we recover the traditional (i.e., the absolute) GIT.
\end{rem}

\begin{say}
The structure of the quotient of $Y^{ss}(\pi^* L^n \otimes M)$ %/\!/G'$
 by the large group $G'$ may be divided
in two ways: mod out by the group $G_0$ first and then by the group $G$; or the
other way
around. In the end, we would have to prove that the quotient is independent of
the order of the
two procedures. But we shall take a unified approach as follows.
\end{say}

\begin{thm}
\label{thm:generalrgit} Keep the assumption as in Theorem
\ref{thm:generalcomparisonofstabilities}.
Then the categorical quotient
$Y^{ss}(\pi^* L^n \otimes M)/\!/G'$ exists
and there is a naturally induced morphism
$$\hat{\pi}: Y^{ss}(\pi^* L^n \otimes M)/\!/G' \rightarrow X^s(L)/\!/G.$$
Moreover,
given an orbit $G \cdot x$ in $X^s(L)$,
the fiber $\hat{\pi}^{-1}([ G \cdot x])$ can be identifibed with
the quotient $((Y_x)^{ss}_{G_0}(M|_{Y_x})/\!/G_0)G_x =
(Y_x)^{ss}_{G_0}(M|_{Y_x})/\!/G_0G_x$ where $G_x$ is the (finite)
isotropy subgroup at the point $x$.
\end{thm}
\begin{pf}  The idea of proof is similar to that of Theorem
\ref{thm:inducedmapfromG-equivariancy}.
We shall present the detail as follows.

By \cite{GIT}, the categorical quotient
$Y^{ss}(\pi^* L^n \otimes M)/\!/G'$ exists as a quasi-projective variety.
By Theorem \ref{thm:generalcomparisonofstabilities} (1), the morphism $\pi$
restricts
to a morphism (still denoted by $\pi$)
$$\pi: Y^{ss}(\pi^* L^n \otimes M) \rightarrow X^{ss}(L).$$
Thus we get an obvious induced morphism
$$\pi: Y^{ss}(\pi^* L^n \otimes M) \rightarrow X^{ss}(L) \rightarrow
X^{ss}(L)/\!/G.$$
Since a categorical quotient is a universal quotient (\cite{GIT}),
the above morphism passes to the quotient to give us the desired morphism
$$\hat{\pi}: Y^{ss}(\pi^* L^n \otimes M)/\!/G' \rightarrow X^s(L)/\!/G.$$

To prove the rest of the statement,
let $[G' \cdot y_1]$ and $[G' \cdot y_2]$ be arbitrary two points in
$Y^{ss}(\pi^* L^n \otimes M)/\!/G'$ that are mapped down to $[G \cdot x] \in
X^s (L)/\!/G$. Then we must have that $\pi (y_1), \pi (y_2) \in X^s (L)$.
Hence $G \cdot \pi(y_1) = G \cdot \pi(y_2) = G \cdot x$. Thus we may well
assume that
$y_1, y_2 \in \pi^{-1}(x)$ by applying some elements of $G$ to $y_1$ and $y_2$
if necessary.
But, in this case,  $[G' \cdot y_1]=[G' \cdot y_2]$
if and only if $G_x [G_0 \cdot y_1]=G_x [G_0 \cdot y_2]$ where $G_x$ acts
naturally
on $(Y_x)^{ss}_{G_0}(M_x)/\!/G_0$ and $[G_0 \cdot y_1]$, $[G_0 \cdot y_2]$ are
regarded as points
in $(Y_x)^{ss}_{G_0}(M_x)/\!/G_0$.  This implies that
 $\hat{\pi}^{-1}([G \cdot x]) \cong (Y_x)^{ss}_{G_0}(M_x)/\!/G_0G_x$.
\end{pf}

We shall need the following technical theorem in our later study of universal
moduli problems.

\begin{thm}
\label{thm:Gmainthmforrelativemoduli}
Let $\pi:Y \rightarrow X$ be a projective morphism between two
algebraic varieties and  equivariant with respect to $G$-actions on both $X$
and $Y$.
Assume that  with respect to some $G$-equivariant projective embedding, $X$ is
\begin{enumerate}
\item  contained in the stable locus; and
\item is closed in the semistable locus.
\end{enumerate}
Then there exists in addition a
$G$-equivariant relative
projective embedding for $\pi:Y \rightarrow X$ and an induced $G$-linearization
such that  $Y$ is
\begin{enumerate}
\item  contained in the stable locus; and
\item  is closed in the semistable locus.  Moreover
\item  $Y/G$ factors naturally to $X/G$ with the
fiber over $[G \cdot x] \in X/G$ isomorphic to $\pi^{-1}(x)/G_x$ for $x \in X$.
\end{enumerate}
\end{thm}
\begin{pf}  It follows from Theorem \ref{thm:inducedmapfromG-equivariancy} by
restricting everything to $\pi:Y \rightarrow X$.
\end{pf}

A  more general technical theorem which we shall quote in studying moduli
problems later
is as follows:

\begin{thm}
\label{thm:G'toGmainthmforrelativemoduli}
Let $\pi:Y \rightarrow X$ be a projective morphism between two
algebraic varieties and  equivariant with respect to a homomorphism $\rho:G'
\rightarrow G$
where $G'$ acts $Y$ and $G$ acts on $X$.
Assume that  with respect to some $G$-equivariant projective embedding, $X$ is
\begin{enumerate}
\item  contained in the stable locus; and
\item is closed in the semistable locus.
\end{enumerate}
Then there exists in addition a
$\rho$-equivariant relative
projective embedding for $\pi:Y \rightarrow X$ and an induced
$G'$-linearization
such that  $Y$ is
\begin{enumerate}
\item  contained and closed in the $G_0$-semistable locus;
\item  contained and closed in the $G$-stable locus;
\item  contained and closed in the $G'$-semistable locus;
\item  $Y/\!/G'$ factors naturally to $X/G$ with the fiber over $[G \cdot x]
\in X/G$
        isomorphic to $\pi^{-1}(x)/\!/G_0G_x$
       for $x \in X$.
\end{enumerate}
\end{thm}
\begin{pf}  It follows from Theorem \ref{thm:generalrgit} by restricting
everything
to $\pi:Y \rightarrow X$.
\end{pf}

\begin{rem} In studying moduli problems, the above two theorems will allow us
to avoid actually writing
down a relative projective embedding (for a universal curve or a relative Quot
scheme) which
is sometimes technical and time consuming.
\end{rem}

\begin{say}
To close this section, let's look at an interesting special case.
Assume that $G$ acts on two projective varieties $X$ and $X_0$ repectively
and $G_0$ acts on $X_0$ only (acts on $X$ trivially).
Set $Y=X_0 \times X$ acted on by $G'= G_0 \times G$
where $G_0$ operates on the first factor only, while $G$ acts diagonally.
Thus the projection onto the first factor $\pi: Y \rightarrow X$ is equivariant
with respect to the projection onto the  second factor $G' \to G$.
(This example has been worked out by R. Pandharipande in the case that
$X_0 = {\bf P}(V)$ and $X= {\bf P}(W)$.)
We begin with considering the diagonal action of $G$ on $Y=X_0 \times X$.
The GIT problem for $G_0$ on Y is the same as the GIT problem for $G_0$ on
$X_0$.

We have the following inclusion
$$\text{Pic}^G(X) \otimes \text{Pic}^G(X_0) \subset \text{Pic}^G(X_0 \times
X).$$
We shall concentrate on the linearizations coming from $\text{Pic}^G(X) \otimes
\text{Pic}^G(X_0)$.
Take $L \in \text{Pic}^G(X)$ and $L_0 \in \text{Pic}^G(X_0)$.
By Theorem \ref{thm:generalizedReichstein}, we have that there exists $m_0$
such that
when $m \ge m_0$,
\begin{enumerate}
\item $Y^{ss}_G(L_0 \otimes L^{\otimes m}) \subset \pi^{-1}(X^{ss}(L))$.
\item $Y^s_G(L_0 \otimes L^{\otimes m}) \supset \pi^{-1}(X^s(L))$.
\end{enumerate}

In case that $X^{ss}(L))=X^s(L)$, we have
$$Y^{ss}_G(L_0 \otimes L^{\otimes m}) = Y^s_G(L_0 \otimes L^{\otimes m})
= \pi^{-1}(X^s(L)) = X_0 \times X^s(L).$$
Clearly, this implies that when $L_0 \in \text{Pic}^{G'}(X_0)$
$$Y^{ss}(L_0 \otimes L^{\otimes m}) = (X_0)_{G_0}^{ss}(L_0)  \times X^s(L).$$
So we get the induced morphism
$$\hat{\pi}: Y^{ss}(L_0 \otimes L^{\otimes m})/\!/G' \rightarrow X^s(L)/\!/G$$
with the fiber isomorphic to  $(X_0)_{G_0}^{ss}(L_0)/\!/G_0G_x$
at the point $[G\cdot x] \in  X^s(L)/\!/G$. This induced morphism $\hat{\pi}$
needs not
to be a trivial fibration. %even if $G$ acts freely on $X^s(L)$.
\end{say}

%\begin{prop} L
%\end{prop}

\section{Generalized Kempf-Ness's theorem}

\begin{say}
Again, we treat $\pi: Y \to X$ as a $G$-equivariant morphism alone.
Place ourselves in the situation of  \ref{say:setupphieta}.
According to Theorem \ref{thm:comparingstabilities},
in the case that  $X^{sss}(\widetilde{\omega}) = \emptyset$,
$\widetilde{\omega}$ determines uniquely a semistable locus upstairs.
This motivates:
\end{say}

\begin{defn} Let $\widetilde{\eta}_0$ be a polarization on the boundary of
${\frak E}^G(Y)$
which is the pullback of $\widetilde{\omega}$ by the morphism $\pi$. Assume
that
$X^{sss}(\widetilde{\omega}) = \emptyset$. Then we set
$Y^{ss}(\widetilde{\eta}_0) = Y^s(\widetilde{\eta}_0)=
\pi^{-1}(X^{ss}(\widetilde{\omega}))$.
\end{defn}

\begin{lem}
\label{lem:forkempf-nessthm}
If $y \in (\Phi^{\widetilde{\eta_0}})^{-1} (0)$,
then $G \cdot y \cap (\Phi^{\widetilde{\eta_0}})^{-1} (0) = K \cdot y$.
\end{lem}
\begin{pf}
First observe that $y \in (\Phi^{\widetilde{\eta_0}})^{-1} (0)$ implies
$\Phi^{\widetilde{\omega}} (\pi (y)) = 0$. That is, $\pi (y) \in
(\Phi^{\widetilde{\omega}})^{-1}(0)$.

We shall adopt the proof of Lemma 7.2 of \cite{Kirwan84}.
Suppose that $g \in G$ is such that $g \cdot y \in
(\Phi^{\widetilde{\eta}_0})^{-1} (0)$.
We want to show that there exists an element $k \in K$ such that $g \cdot y = k
\cdot y$.
Since $\Phi^{\eta_0}$ is $K$-equivariant and $G=K \text{exp}(i {\frak k})$,
we may assume that $g = \text{exp}(i a), a \in {\frak k}$.
%Thus we obtain $g \cdot \pi(y) \in (\Phi^{\widetilde{\omega}})^{-1}(0)$.
%But since  $X^{sss}(\widetilde{\omega}) = \emptyset$, we already have
%$G \cdot \pi(y) \cap (\Phi^{\widetilde{\omega}})^{-1}(0) = K \pi (y)$ (see
%%\cite{Ness}
or Lemma 7.2 of
Set $$h(t) = \Phi^{\widetilde{\eta}_0} (\text{exp}(i at) \cdot y) \cdot a.$$
Then $$h(t) = \Phi^{\widetilde{\omega}} \circ \pi (\text{exp}(i at) \cdot y)
\cdot a.$$
$h(t)$ vanishes at $t=0$ and $t=1$ because both $y$ and $\text{exp}(i a) \cdot
y$
belong to  $(\Phi^{\widetilde{\eta}_0})^{-1} (0)$. Thus there must be a point
 $t \in (0, 1)$ such that
$$0 = h'(t) = d \Phi^{\widetilde{\omega}} \circ d \pi (i (\xi^Y_a)_z) \cdot a
= \omega (d \pi (i (\xi^Y_a)_z), (\xi^X_a)_{\pi(z)})$$
$$= \omega (i (\xi^X_a)_{\pi(z)}, (\xi^X_a)_{\pi(z)}) =<(\xi^X_a)_{\pi(z)},
(\xi^X_a)_{\pi(z)}>$$
where $z = \text{exp}(i at) \cdot y$.
Hence $(\xi^X_a)_{\pi(z)} = 0$. Or, $a \in \text{Lie}(G_{\pi(z)})$. But
$$\pi(z) = \text{exp}(i at) \cdot \pi (y) \in G
(\Phi^{\widetilde{\omega}})^{-1}(0).$$
So $a=0$ because $G_{\pi(z)}$ must be a finite group. This completes the proof.
\end{pf}

Now we have the following theorem which extends the Kempf-Ness theorem to the
degenerated polarizations.

\begin{thm} {\rm \text{(Generalized Kempf-Ness's Theorem)}}
Assume that $X^{ss}(\widetilde{\omega})=X^s(\widetilde{\omega})$. Then
\begin{enumerate}
\item $Y^{ss}(\widetilde{\eta}_0) = G (\Phi^{\widetilde{\eta}_0})^{-1} (0)$;
\item the topological quotient $Y^{ss}(\widetilde{\eta}_0)/G$ is Hausdorff;
\item the inclusion $(\Phi^{\widetilde{\eta}_0})^{-1} (0)
\hookrightarrow Y^{ss}(\widetilde{\eta}_0)$ induces
a homeomorphism between \\
$(\Phi^{\widetilde{\eta}_0})^{-1} (0)/K$ and $Y^{ss}(\widetilde{\eta}_0)/G$.
\item $Y^{ss}(\widetilde{\eta_0})/G = (\Phi^{\widetilde{\eta}_0})^{-1} (0)/K$
inherits from $\eta_0$
a closed 2-form (which may be degenerated somewhere) away from singularities.
\end{enumerate}
\end{thm}

\begin{pf}  (1). By definition,
$$Y^{ss}(\widetilde{\eta}_0) = \pi^{-1}(X^s(\widetilde{\omega})
= \pi^{-1}(G (\Phi^{\widetilde{\omega}})^{-1}(0)  = G
\pi^{-1} (\Phi^{\widetilde{\omega}})^{-1}(0) $$
$$= G (\Phi^{\widetilde{\omega}} \circ \pi)^{-1}(0)
= G (\Phi^{\widetilde{\eta_0}})^{-1} (0).$$

(2). This follows from the fact that $Y^{ss}(\widetilde{\eta}_0)=
Y^{ss}(\widetilde{\eta}(t))$
for sufficiently small positive numbers $t$.

(3). (1) implies that the induced map $(\Phi^{\eta_0})^{-1} (0)/K \to
Y^{ss}(\widetilde{\eta}_0)/G$
is surjective. Lemma \ref{lem:forkempf-nessthm} implies that it is injective.
Now as a continuous bijection between Hausdorff spaces, it must be a
homeomorphism.

(4) The proof is the same as the one for non-degenerated 2-forms
(\cite{MarsdenWeinstein}).
\end{pf}

\begin{rem} In fact by Theorem \ref{thm:comparingstabilities},
  the quotient $Y^{ss}(\widetilde{\eta}_0)/G$
admits a complex structure and many other non-degenerated 2-forms induced from
$\widetilde{\eta}(t), t \in (0, \delta]$ such that the 2-form induced from
$\widetilde{\eta}_0$
is the limit of the above.
\end{rem}

\begin{rem}
When $X^{sss}(\widetilde{\omega}) \ne \emptyset$, there is a difficulty in
defining that
$$Y^{ss}(\widetilde{\eta}_0) := \{ y \in Y | 0 \in \Phi^{\widetilde{\eta}_0}
(\overline{G \cdot y}) \}
= \{ y \in Y | 0 \in \Phi^{\widetilde{\omega}} (\overline{G \cdot \pi (y)}) \}
= \pi^{-1}(X^{ss}(\widetilde{\omega})).$$
The problem is as follows. Let $x \in X^{sss}(\widetilde{\omega})$ be a point
such that
$G \cdot x$ is closed in $X^{ss}(\widetilde{\omega})$. Then $G_x$ is reductive
and acts
on the fiber $\pi^{-1}(x) \subset Y^{ss}(\widetilde{\eta}_0)$. Since $\dim G_x
> 0$,
$\pi^{-1}(x)/ G_x$ is non-Hausdorff. Thus one would like to exclude some points
in
$\pi^{-1}(x)$ from $ Y^{ss}(\widetilde{\eta}_0)$, presumly by using a
non-degenerate
closed 2-form on $\pi^{-1}(x)$.  That is, one would like to pick up some
semistability on  $\pi^{-1}(x)$ for the action of $G_x$.
The form $\eta_0$ is helpless in this regard
since the fiber $\pi^{-1}(x)$ is exactly where $\eta_0$ vanishes.
Hence among many choices of the semistabilities on  $\pi^{-1}(x)$ for the
action of $G_x$,
we do not know, a priori, which to choose.
Notice that the same ambuguity does not happen when
$X^{sss}(\widetilde{\omega}) = \emptyset$.
In this case, the isotropy subgroup $G_x$ for every $x \in
X^{ss}(\widetilde{\omega})$
is a finite group. $\pi^{-1}(x)/ G_x$ is always a good quotient.

One of ways to solve (actually to pass) the problem of
$X^{sss}(\widetilde{\omega}) \ne \emptyset$ is as follows. Assume that $\dim
{\frak E}^G(X) > 1$
and it has at least one top chamber. Choose a linearization
$\widetilde{\omega}'$
in a top chamber that contains $\widetilde{\omega}$ in its closure.
We then have $X^{ss}(\widetilde{\omega}') = X^s(\widetilde{\omega}') \subset
X^{ss}(\widetilde{\omega})$.
Now $\pi^{-1}(X^{ss}(\widetilde{\omega}'))$ has a good complete quotient
and one has natural maps
$$\pi^{-1}(X^{ss}(\widetilde{\omega}'))/\!/G \rightarrow
X^{ss}(\widetilde{\omega}')/\!/G \rightarrow
X^{ss}(\widetilde{\omega})/\!/G .$$
This helps to reduce the case when $X^{sss}(\widetilde{\omega})\ne \emptyset$
to the nicer case when  $X^{sss}(\widetilde{\omega}) = \emptyset$.
\end{rem}

\begin{rem}
In general, the semistable set $Y^{ss}(\pi^*L^n \otimes M) (n \gg 0)$ may be
recovered as follows:
Let $X^{ss}_c (L)$ be the set of closed orbits in $X^{ss}(L)$. Then
\begin{enumerate}
\item  $Y_c^{ss}(\pi^*L^n \otimes M) = \{ y \in Y | y \in \pi^{-1}
(X_c^{ss}(L)) \cap
(Y_{\pi (y)})_c^{ss} (M_{\pi (y)}) \}$.
\item  $Y^{ss}(\pi^*L^n \otimes M) =
\{ y \in Y | \overline{G \cdot y} \cap Y_c^{ss}(\pi^*L^n \otimes M) \ne
\emptyset \}$.
\item  $Y^s (\pi^*L^n \otimes M)= \{y \in Y | y \in Y_c^{ss}(\pi^*L^n \otimes
M),
G_y \; \text{ is finite} \}$.
\end{enumerate}
\end{rem}

\begin{rem}
There are nef
linearized line bundles that are not  the pull-backs for any algebraic
contraction maps.
In this case, we do not know yet how to make sense of RGIT in the algebraic
category.
One possible alternative is to allow the contractions
to be just complex analytic and work in the complex analytic category. The
price paid
is then the loss of algebraicity and the possible validity over positive
characteristics.
\end{rem}

\section{Configuration spaces and Grassmannians}

In this section, we shall present an application of our various RGIT theorems
to some
elementary finitely dimensional settings (as opposed to the later applications
to
the moduli problems).

\begin{say}
Consider the Grassmannian of $n$-subspaces in ${\Bbb C}^m$,
$\text{Gr}(n, {\Bbb C}^m)$, acted on by the maximal torus $T=({\Bbb
C}^*)^{m-1}$.
Since $\text{Pic Gr}(n, {\Bbb C}^m) \cong {\Bbb Z}$ and is generated by an
ample line bundle,
the stabilities of
linearized actions are determined by the characters of  $T=({\Bbb C}^*)^{m-1}$.
In fact, let $\Phi: \text{Gr}(n, {\Bbb C}^m) \rightarrow {\Bbb R}^m$ be the
standard moment map induced by the Pl\"uker embedding.
Then the moment map image is the so called hypersimplex $\Delta^m_n$
$$\Delta^m_n = \{(\alpha_1, \cdots \alpha_m) | 0 \le \alpha_i \le 1, \sum_j
\alpha_j = n\}.$$
Thus we have that the $T$-effective ample cone
 ${\frak E}^T(\text{Gr}(n, {\Bbb C}^m))$
can be identified with the cone over  $\Phi(\text{Gr}(n, {\Bbb
C}^m))=\Delta^m_n$.
\end{say}

\begin{say}
On the other hand, consider the diagonal action of $G=\text{SL}(n)$ on
$X=({\Bbb P}^{n-1})^m$
$$\text{SL}(n) \times ({\Bbb P}^{n-1})^m \rightarrow ({\Bbb P}^{n-1})^m.$$
We have $\text{Pic}^G(X)
\cong \text{Pic}(X) \cong {\Bbb Z}^m$. A nef line bundle
$$L \cong {\cal O}_{{\Bbb P}^{n-1}}(k_1)
\otimes \cdots \otimes {\cal O}_{{\Bbb P}^{n-1}}(k_m), \;\; k_i \ge 0 $$
is $G$-effective (i.e., $X^{ss}(L) \ne \emptyset$)
 if and only if $n k_i \leq \sum_j k_j$
for any $i = 1, \cdots, m$. If we set $\alpha_i = n k_i/\sum_j k_j$,
we can express this condition by
the inequalities $ 0 \le \alpha_i \leq 1, \sum \alpha_i = n$.
This shows that the $G$-effective ample cone ${\frak E}^G(X)$
 is equivalent  to the cone over the polytope
$$\Delta^m_n = \{ (\alpha_1, \cdots, \alpha_m) \in
{\Bbb R}^m: 0 \le \alpha_i \leq 1, \sum \alpha_i = n\}.$$
\end{say}

\begin{say}
\label{say:GMcorrespondence}
Hence we obtain the identification ${\frak E}^T(\text{Gr}(n, {\Bbb C}^m)) =
{\frak E}^G(({\Bbb P}^{n-1})^m)$.
Applying the 1-1 correspondence between the $T$-orbits on $\text{Gr}(n, {\Bbb
C}^m)$
and the $G$-orbits on $({\Bbb P}^{n-1})^m$ (\cite{GelfandMacPherson}),
we can further obtain the identification
bewteen the $T$-GIT quotients on the Grassmannian $\text{Gr}(n, {\Bbb C}^m)$
and the $G$-GIT quotients on $({\Bbb P}^{n-1})^m$.
\end{say}

\begin{say}
Notice that the underlying line bundle of any element in ${\frak
E}^T(\text{Gr}(n, {\Bbb C}^m))$
is ample. Thus the boundary of ${\frak E}^T(\text{Gr}(n, {\Bbb C}^m))$ consists
of
degenerating characters but with ample underlying line bundles.
However, on the contrast,  notice that the group $\text{SL}(n)$ has no
characters and the
boundary of ${\frak E}^G(X)$ consists of only {\sl nef} line bundles.
This shows an interesting phenomenon that
the two sorts of degenerations can sometimes be harmoniously linked.
\end{say}

\begin{say} In the following,
we shall use the above identifications freely in the rest of paper. One will
see that
one point of view sometimes has advantage over the other (and vice versa).
We will use the above to exhibit our theorems in \S\S 3 and 4.
To simplify exhibition, we only consider the case when $n=2$, i.e.,
the action of $T=({\Bbb C}^*)^{m-1}$ on $\text{Gr}(2, {\Bbb C}^m)$ and
the action of $G=\text{SL}(2)$ on $({\Bbb P}^1)^m$, leaving out some possible
generalizations for $n > 2$ to a future paper.
\end{say}

\begin{say}
{}From \ref{say:GMcorrespondence}, we have
$${\frak E}^T(\text{Gr}(2, {\Bbb C}^m)) = {\frak E}^G(({\Bbb P}^1)^m)
=\Delta^m_2.$$
The walls of $\Delta^m_2$ are of the form
$$W_J = \{ (\alpha_1, \cdots, \alpha_m) \in \Delta^m_2 | \sum_{i \in J}
\alpha_i = 1\}$$
where $J \subset \{1, \cdots, m\}$ is any proper subset.
The faces of $\Delta^m_2$  but vertices
can be obtained by setting some of the coordinates to be 0 or
one coordinate to be 1. (If two coordinates are 1, we will get a vertex.)
Thus they are divided into two types: the ones that are obtained by
setting $k$ coordinates to be 0 are again hypersimplexes $\Delta^{m-k}_2$;
the ones that are obtained by
setting $k-1$ coordinates to be 0 but one coordinate to be 1 become simplexes
$\Delta^{m-k}_1$.

In particular, there are in general two different types of facets (faces of
codimension 1):
$\Delta^{m-1}_2$ and $\Delta^{m-1}_1$. Precisely, they are
$$\Delta^{m-1}_2 [i] =
\{ (\alpha_1, \cdots, \alpha_m) \in \Delta^m_2 | \alpha_i = 0, \sum_{j \ne i}
\alpha_j = 2\},
\; \text{for all} \; 1 \le i \le m.$$
$$\Delta^{m-1}_1 [i] =
\{ (\alpha_1, \cdots, \alpha_m) \in \Delta^m_2 |  \alpha_i = 1, \sum_{j \ne i}
\alpha_j = 1\},
\; \text{for all} \; 1 \le i \le m.$$
\end{say}

\begin{say}
 Given an element $\alpha=  (\alpha_1, \cdots, \alpha_m) \in \Delta^m_2$,
we use ${\cal M}^m_\alpha$ to denote the corresponding (GIT or symplectic)
quotient
of either  $({\Bbb P}^1)^m$ by $G=\text{SL}(2)$ or $\text{Gr}(2, {\Bbb C}^m)$
by $T=({\Bbb C}^*)^{m-1}$.
\end{say}

\begin{rem}
\label{say:polygons}
We point out that the (GIT or symplectic) quotients
of these two actions admit some other interesting interpretations. First of
all,
they appear as the arithmetic quotients of complex balls (Deligne-Mostow).
Secondly,
they are also the moduli spaces of spatial polygons modulo
orientation-preserving Euclidean
motions (\cite{Klyachko92}, \cite{KapovichMillson94b}). Moreover, it is easy to
see that
the real points of these quotients are the trivial double covers (disjoint
unions) of the
moduli spaces of linkages in Euclidean palne  modulo orientation-preserving
Euclidean
motions (\cite{KapovichMillson94a} and \cite{Hu92}). Thus,
 as the sets of the real points of the symplectic quotients ${\cal
M}^m_\alpha$,
the topological changes of the moduli spaces of linkages with prescribed side
lengths
in the Euclidean plane when crossing walls are governed by the changes of
${\cal M}^m_\alpha$
(cf. \cite{Hu92}). This is to say that two moduli spaces of linkages
(without zero side lengths) in the Euclidean plane are related by a sequence of
real
blowup and blowdowns.
\end{rem}

\begin{say}
There are $m$ many forgetful morphisms from $({\Bbb P}^1)^m$ to
$({\Bbb P}^1)^{m-1}$. Use $f_i$ to denote the  forgetful map obtained
by forgetting the $i$th factor from  $({\Bbb P}^1)^m$. These are
$G$-equivariant {\sl trivial}
fibrations
with fiber ${\Bbb P}^1$.   In terms of $G$-effective cones, $f_i$
corresponds to the inclusion of the facet
$$\Delta^{m-1}_2 [i] \subset \Delta^m_2$$
where $\Delta^{m-1}_2 [i] \cong \Delta^{m-1}_2$
is obtained by setting the $i$th coordinate to be zero.
Noting that ${\cal M}^m_{(\alpha_1, \cdots,\alpha_{i-1}, 0,
\alpha_{i+1}, \cdots, \alpha_m)} = {\cal M}^{m-1}_{(\alpha_1,
\cdots,\alpha_{i-1},
\alpha_{i+1}, \cdots, \alpha_m)}$, we have
\end{say}

\begin{thm} (cf. Theorem 2.1, \cite{Hu92})
\label{thm:Hu92thm2.1}
 Let $\alpha= (\alpha_1, \cdots,\alpha_{i-1},
\alpha_{i+1}, \cdots, \alpha_m) \in \Delta^{m-1}_2$
do not lie on any wall and set
$$\widetilde{\alpha}_\epsilon = (\alpha_1 -{\epsilon \over {m-1}}, \cdots,
\alpha_{i-1}-{\epsilon \over {m-1}},
\epsilon, \alpha_{i+1}-{\epsilon \over {m-1}}, \cdots,  \alpha_m-{\epsilon
\over {m-1}})
 \in \Delta^m_2$$
where $\epsilon$ is a sufficiently small
 positive number. Then
${\cal M}^m_{\widetilde{\alpha}_\epsilon} = {\cal M}^{m-1}_\alpha \times {\Bbb
P}^1$.
\end{thm}
\begin{pf}
 Applying Theorem \ref{thm:inducedmapfromG-equivariancy} (2) to the forgetful
map
$f_i: ({\Bbb P}^1)^m \rightarrow ({\Bbb P}^1)^{m-1}$,
one sees immediately that ${\cal M}^m_{\widetilde{\alpha}_\epsilon}$ is a
fibration over
${\cal M}^{m-1}_\alpha$ with fibers ${\Bbb P}^1$. The triviality follows from
the fact that $f_i$ is equivariantly trivial and $\text{SL}(2)$ acts freely
on any stable locus.
\end{pf}

\begin{rem}
Using the remark in \ref{say:polygons} about the real loci of the quotients
${\cal M}^m_\alpha$, Theorem \ref{thm:Hu92thm2.1}
 gives an alternative justification for Corollary 15 of
\cite{KapovichMillson94a}.
\end{rem}

\begin{rem}
In terms of Grassmannians, the  forgetting map $f_i$
corresponds to the (rational) projection  from $\text{Gr}(2, {\Bbb C}^m)$ to
 $\text{Gr}(2, {\Bbb C}^{m-1})$ by projecting a 2-plane in ${\Bbb C}^m$ to
the $i$th coordinate $(m-1)$-hyperplane (the hyperplane is
obtained by setting the $i$th coordinate to be zero).
\end{rem}

\begin{say}
To study the inclusion $$\Delta^{m-1}_1 [i] \subset \Delta^{m}_2$$
where $\Delta^{m-1}_1 [i]$ is obtained by setting the $i$th coordinate to be 1,
we switch our point of view from configuation spaces of points on  ${\Bbb P}^1$
to Grassmannians.
There are $m$ many {\sl rational} facet maps ${\frak f}_i$
from $\text{Gr}(2, {\Bbb C}^m)$ to
 $\text{Gr}(1, {\Bbb C}^{m-1})={\Bbb P}^{m-2}$ by taking the intersection of
a  2-plane in ${\Bbb C}^{m}$ with the $i$th coordinate $(m-1)$-hyperplane
(the hyperplane is obtained by setting the $i$th coordinate to be zero).
These are {\sl truly} rational maps but are equivariant with respect to the
projection
$\rho$ from the maximal torus
$T= ({\Bbb C}^*)^{m-1} = ({\Bbb C}^*)^m/\hbox{(diagonal)}$ to the quotient
group $T/T_1$
where the 1-PS $T_1 = \{(0, \cdots, \lambda,\cdots, 0) | \lambda \in {\Bbb C}^*
\}$ (modulo the
diagonal group) acts trivially on  $\text{Gr}(1, {\Bbb C}^{m-1})={\Bbb
P}^{m-2}$.
Any element of $\Delta^{m-1}_1 [i]$ is of the form
$\alpha = (\alpha_1, \cdots,\alpha_{i-1},
1, \alpha_{i+1}, \cdots,  \alpha_m)$. Since $T/T_1= ({\Bbb C}^*)^{m-2}$ acts on
$\text{Gr}(1, {\Bbb C}^{m-1})={\Bbb P}^{m-2}$ with a dense open orbit,
we see that ${\cal M}^m_{\alpha}$ is a point for such an element $\alpha$.
\end{say}

\begin{thm} (cf. Theorem 2.1, \cite{Hu92}) Let $\alpha = (\alpha_1,
\cdots,\alpha_{i-1},
1, \alpha_{i+1}, \cdots,  \alpha_m)$ be any element in the interior of
 $\Delta^{m-1}_1 [i]$ and set
$$\widetilde{\alpha}_{1-\epsilon} = (\alpha_1 +{\epsilon \over {m-1}},
 \cdots,\alpha_{i-1}+{\epsilon \over {m-1}},
1- \epsilon, \alpha_{i+1}+{\epsilon \over {m-1}}, \cdots,
\alpha_{m+1}+{\epsilon \over {m-1}})$$ where $\epsilon$ is a sufficiently small
 positive number.
Then  ${\cal M}_{\widetilde{\alpha}_{1-\epsilon}}$ is isomorphic to ${\Bbb
P}^{m-3}$.
\end{thm}

\begin{pf} By Bialynicki-Birula's decomposition theorem
\cite{Bialynicki-Birula73},
the fibers of the map ${\frak f}_i$ have the structures of $T_1$-modules where
$T_1$ acts with positive weights. Since $T$ acts quasi-freely (no finite
isotropy subgroups
except for the identity group), all of these weights are 1.
The theorem now follows from  Theorem \ref{thm:generalrgit}
 and the fact that
 ${\cal M}^m_{\alpha}$ is a point.
\end{pf}

\begin{rem} In terms of configurations of points on ${\Bbb P}^1$, the facet map
${\frak f}_i$ correspond to the projection onto the $i$th factor ($1 \le i \le
m$).
The reason that we hesitated working over $({\Bbb P}^1)^m$ is that
the  $i$th factor (${\Bbb P}^1$) has no stable points with respect to the
action
of $\hbox{SL}(2)$ although a quotient trivially exists.
\end{rem}

%\smallskip
%\noindent
%3.1.6 Let $G$ be a reductive algebraic group and $B$ a fixed Borel subgroup.
%Then $Y=G/B$ is the largest compact homogeneous variety that factors to any
%%other
%compact homogeneous variety whose automorphism group is $G$. In this section,
%we shall focus our attentions on the degenerated GIT on $G/B$ and its relation
%with the relative GIT for the {\sl regular} equivariant morphisms from
%$G/B$ to some other  homogeneous variety $G/P$, where $P$ is a parabolic
%%subgroup
%containing the Borel subgroup $B$.

\vfill\eject

\begin{center}
{\bf II. Universal Moduli}
\end{center}

%\section{Review of the base moduli spaces  $\overline{M_{g}}$ and
% $\overline{M_{g, n}}$}
\section{Quest for universal moduli spaces}

\begin{say}
  \label{say:moduli spaces for curves}
  Recall that one has the following coarse moduli spaces in the case of curves.
  \begin{enumerate}
  \item $M_{g}$, parameterizing nonsingular curves of genus $g\ge2$
    and its compactification $\overline{M_{g}}$, parameterizing
    Mumford-Deligne stable curves. %see Mumford \cite{Mumford77},
  \item  $M_{g,m}$, $2g-2+m>0$, parameterizing nonsingular
$m$-pointed curves of genus $g\ge2$,
 and its compactification $\overline{M_{g,m}}$,
 parameterizing  stable   $m$-pointed curves. %see Knudsen \cite{Knudsen83},
  %\item a moduli space $M_{g,m}(W)$ of stable maps from reduced curves
  %  to a variety $W$, see Kontsevich \cite{Kontsevich94}.
  \end{enumerate}
  It is well known that $\overline{M_{g}}$ and $\overline{M_{g,m}}$ are
  projective.
\end{say}

Precisely,
the Mumford-Deligne stable curves (resp. $m$-pointed curves) are defined as
follows:
\begin{defn}
\label{defn:definitionofstablecurves}
Let $S$ be a base scheme.
 A stable (resp. semistable) curve of genus $g \ge 2$ over $S$ is a proper flat
morphism
$f: C \to S$ such that for all $s \in S$ the geometric fiber $C_s$ of $f$ over
$s$,
satisfies:
\begin{enumerate}
\item $C_s$ is reduced, connected scheme of dimension 1 with $h^1(C_s, {\cal
O}_{C_s}) =g$;
\item every singular point of $C_s$ is an ordinary double point;
\item if $C_1$ is an irreducible rational component of $C_s$ then $C_1$ meets
the rest of
$C_s$ in at least 3 points (resp. 2 points).
\end{enumerate}
\end{defn}

\begin{defn}
\label{defn:definitionofstablen-pointedcurves}
 A $m$-pointed curve stable (resp. semistable)
curve of genus $g \ge 2$ over $S$ is a proper flat morphism
$f: C \to S$ together with $m$-distinct sections $s_i:S \to C$ such that
for all $s \in S$ the geometric fiber $C_s$ of $f$ over $s$,
satisfies:
\begin{enumerate}
\item $C_s$ is reduced, connected scheme of dimension 1 with $h^1(C_s, {\cal
O}_{C_s}) =g$;
\item every singular point of $C_s$ is an ordinary double point;
\item all $s_i(s)$ are smooth points of $C_s$ and $s_i(s)\ne s_j(s)$ for $i \ne
j$
\item if $C_1$ is an irreducible rational component of $C$ then the number of
points
where $C_1$ meets the rest of
$C_s$ plus the number of points  $s_i(s)$ on $C_1$ is at least 3 (resp. 2).
\end{enumerate}
\end{defn}

Rather than going through some abstract definitions of universal moduli spaces,
let us go over (briefly) the concrete universal moduli problems listed in 0.1
of the introduction.

\begin{say}
\label{say:listofmoduli}
\begin{enumerate}
%\item The universal moduli space $\overline{CM}_{g,n}$ of $n$-Cartesian
%%products of stable curves
%   over the moduli space $\overline{M_g}$ of Delighe-Mumford stable curves
%   of genus $g$. (That is, given a stable
%   curve $[C]$ in  $\overline{M_g}$, the fiber in the universal moduli space
%  over the point $[C]$ is functorially isomorphic to $C^n$ modulo the
%%automorphism
%  group of $C$.)

\item The universal moduli space $\overline{FM_{g, n}} \rightarrow
\overline{M_g}$
      of Fulton-MacPherson configuration spaces of stable curves. (see also
\cite{Pandharipande94b})

\item The compactified
      universal Picard $\overline{P^d_g}  \rightarrow \overline{M_g}$ of degree
$d$ line bundles
      (\cite{Caporaso94}).

\item The universal moduli space $\overline{P_{g,m}(e, r, F, \alpha)}
\rightarrow
     \overline{M_{g,n}}$ of $p$-semistable parabolic sheaves of degree $e$,
rank $r$,
     type $F$,  and weight $\alpha$
      (\cite{Hu95}).

\item The universal moduli space $M_{g}({\cal O}, P) \rightarrow
\overline{M_g}$
 of $p$-semistable coherent sheaves with a fixed Hilbert polynomial $P$
To be more precise, for each $[C] \in \overline{M_{g}}$, one has a natural
projective variety,
$M_C({\cal O}_C, P)$,
parameterizing (Simpson's) $p$-semistable coherent sheaves of a fixed Hilbert
polynomial
$P$. Our aim is to construct the universal moduli space
 $M_g({\cal O}, P) \rightarrow \overline{M_{g}}$ parametrizing
the set of equivalence classes of pairs $(C, E)$ where
$[C] \in \overline{M_{g}}$ and $[E] \in M_C({\cal O}_C, P)$
such that the fiber over the stable curve $C$ is $M_C({\cal O}_C,
P)/\hbox{Aut}(C)$.

\item The universal Hilbert scheme $\hbox{Hilb}^n_g  \rightarrow
\overline{M_g}$
of 0-dimensional schemes of length $n$ on the Mumford-Deligne  stable curves.
There exists a natural
dominating morphism $$\psi: \hbox{Hilb}^n_g  \longrightarrow M_g({\cal O}, P)$$
when $P(x) = x + n + 1 - g$.
\end{enumerate}
\end{say}

\begin{rem} The last two cases in the above list
 are what we shall considered seriously in the sequel.
The construction of the universal moduli space $M_{g}({\cal O}, P) \rightarrow
\overline{M_g}$
 will be done by using our RGIT.
%This space has also been constructed in \cite{Pandharipande94a}. But our
%%approaches
%are much easier and promise  straightforward
%generalizations to the  universal moduli spaces of $p$-semistable coherent
%%sheaves
%over the moduli spaces of higher dimensional varieties (cf. \cite{Vieweg95}),
%%because our construction
%depends on Simpson's construction of moduli space of $p$-semistable coherent
%%sheaves, and his
%construction works for any projective scheme! (To make the material in this
%%paper coherent,
%we shall postpone such generalizations to a later publication. See
%%\cite{Hu95}.)
\end{rem}

\begin{rem} The moduli problem in \ref{say:listofmoduli} (3) can be specified
as follows.
For each $[(C, p_1, \cdots, p_m)] \in M_{g,m}$,
one has a natural projective variety, $P_C(e, r, F, \alpha)$,
parameterizing $\alpha$-semistable vector bundles of degree $e$ and rank $r$
with quasi-parabolic structures of type $F$ at points $p_i$.
Let $P_{g,m}(e, r, F, \alpha)$
be the set of equivalence classes of pairs $((C, p_1, \cdots, p_m),
{\cal E})$ where $[(C, p_1, \cdots, p_m)] \in M_{g,m}$
and $[{\cal E}] \in P_C(e, r, F, \alpha)$. Now it is natural to ask for a
compactification $\overline{P_{g,m}(e, r, F, \alpha)}$ of
$P_{g,m}(e, r, F, \alpha)$ with the following desired properties
(cf. \cite{Pandharipande94a}):
\begin{enumerate}
\item  $\overline{P_{g,m}(e, r, F, \alpha)}$
is a projective variety parameterizing equivalence classes of algebro-geometric
objects.
\item $\overline{P_{g,m}(e, r, F, \alpha)}$
contains $P_{g,m}(e, r, F, \alpha)$ as an open dense subset.
\item There is a natural morphism $\eta: \overline{P_{g,m}(e, r, F, \alpha)}
\rightarrow \overline{M_{g,m}}$ such that
the following natural diagram commutes:
\begin{equation*}
\begin{CD}
P_{g,m}(e, r, F, \alpha) @>>> \overline{P_{g,m}(e, r, F, \alpha)} \\
@VVV @VV{\eta}V \\
M_{g, m} @>>> \overline{M_{g, m}}
\end{CD}
\end{equation*}
\item For each $[(C, p_1, \cdots, p_m)] \in M_{g,m}$, there is an isomorphism
$$\eta^{-1}([(C, p_1, \cdots, p_m)]) \cong P_C(e, r, F, \alpha)/\hbox{Aut}(C,
p_1, \cdots, p_m).$$
\end{enumerate}
Because of the lack of satisfactory GIT construction of $\overline{M_{g,m}}$,
we postpone treating this problem in a later publication \cite{Hu95}.
\end{rem}

\section{The universal moduli space $M_g ({\cal O}, P) \rightarrow
\overline{M_g}$}
\label{section:simpson}
\begin{say} {\sl Simpson's construction of the moduli space of $p$-semistable
coherent sheaves.}
Let $X$ be a projective scheme over ${\Bbb C}$ with a  very ample invertible
sheaf ${\cal O}_X(1)$.
For any coherent sheaf ${\cal E}$ over $X$, let $p({\cal E}, n)$ be the Hilbert
polynomial of
 ${\cal E}$ with $p({\cal E}, n)= \dim H^0(X, {\cal E}(n))$ for $n \gg 0$.
Let $d=d({\cal E})$ be the dimension of the support of ${\cal E}$ which is also
the degree of the Hilbert polynomial of
 ${\cal E}$. The leading coefficient is $r/d!$ where $r=r({\cal E})$ is an
integer which
is called the rank of ${\cal E}$. A coherent sheaf
${\cal E}$ is of pure dimension $d=d({\cal E})$
if for any non-zero subsheaf ${\cal F} \subset
{\cal E}$, we have that $d({\cal E}) = d({\cal F})$.
\end{say}

\begin{defn} A coherent sheaf ${\cal E}$ is $p$-semistable (resp. $p$-stable)
if it is of
pure dimension, and if for any subsheaf ${\cal F} \subset
{\cal E}$, there exists $N$ such that for $n \ge N$
$$\frac{p({\cal F}, n)}{r({\cal F})} \le \frac{p({\cal E}, n)}{r({\cal E})}.$$
\end{defn}

As usual, a $p$-semistable sheaf ${\cal E}$ admits a filtration by subsheaves
$$0 = {\cal E}_0, \subset {\cal E}_1 \subset \cdots \subset {\cal E}_k ={\cal
E}$$
such that the quotient sheaves  ${\cal E}_1/{\cal E}_{i-1}$
are $p$-stable. This filtration is not unique, but
 $gr({\cal E}) = \oplus {\cal E}_1/{\cal E}_{i-1}$ is.
Two  $p$-semistable sheaves ${\cal E}$ and ${\cal E}'$ are
$s$-equivalent if $gr({\cal E}) = gr({\cal E}')$.

Simpson also extends the above definition to the following relative version.

\begin{say} Let $S$ be a base scheme of finite type over ${\Bbb C}$ and $X
\rightarrow S$
 a projective scheme over $S$. Fix a (Hilbert) polynomial $P$ of degree $d$.
A $p$-semistable (resp. stable) sheaf ${\cal E}$ on $X/S$ with Hilbert
polynomial $P$
is a coherent sheaf ${\cal E}$ on $X$, flat over $S$, such that for each
closed point $s \rightarrow S$, ${\cal E}_s$ is a $p$-semistable (resp. stable)
sheaf
of pure dimension $d$ and Hilbert polynomial $P$ on the fiber $X_s$.
\end{say}

\begin{say} {\sl Hilbert schemes and Grassmannians}.
\label{say:hilbertschemes}
As above, let $X \rightarrow S$
be a projective scheme over $S$ with a relatively very ample invertible sheaf
${\cal O}_X(1)$.
Fix a (Hilbert) polynomial $P(n)$. Suppose that ${\cal W}$ is a coherent sheaf
on $X$ flat
over $S$. The Hilbert scheme
$\hbox{Hilb}({\cal W}, P)$ parametrizing quotients
$${\cal W} \rightarrow {\cal F} \rightarrow 0$$
with Hilbert polynomial $P$.  The fiber of $\hbox{Hilb}({\cal W}, P)$ over a
closed point
$s \in S$ is $\hbox{Hilb}({\cal W}_s, P)$.

Grothendieck gives some very explicit
relative projective embeddings of $\hbox{Hilb}({\cal W}, P)$ over
$S$. There is an $M > 0$ such that for any $m \ge M$ we get a closed embedding
$$\psi_m: \hbox{Hilb}({\cal W}, P) \rightarrow \hbox{Grass}(H^0(X/S, {\cal
W}(m)), P(m)).$$
There is a canonical invertible sheaf ${\cal L}_m$
 on the relative Grassmannian by the embedding $\psi_m$.
Over any point in the Grassmannian represented by the quotient ${\cal W}
\rightarrow {\cal F}$,
the restriction of the  invertible sheaf   ${\cal L}_m$ is canonically
identified with
the $\wedge^{P(m)} H^0(X/S, {\cal F}(m))$.
{\sl These invertible sheaves are the ones we shall use in place of the
relatively ample line bundles $M$ in  our various theorems in RGIT} (\S\S 3 and
4).
\end{say}

\begin{say} {\sl The construction of the moduli spaces of $p$-semistable
coherent sheaves}.
\label{say:Q2}
Fix a large number $N$. Let ${\cal W}={\cal O}_X(-N)$ and $V={\Bbb C}^{P(N)}$.
Let $Q_1 \subset \hbox{Hilb}(V \otimes {\cal W}, P)$ denote the open subset of
$p$-semistable
sheaves of pure dimension $d$. We can assume that $N$ is chosen large enough
so that: every $p$-semistable coherent sheaf with Hilbert polynomial $P$
appears as a quotient
corresponding to a point of $Q_1$. Now set $Q_2$ equal to
the open subset in $Q_1$ such that $\alpha: V\otimes {\cal O}_S \rightarrow
H^0(X/S, {\cal E}(N))$ is isomorphism where $\alpha$ is a morphism such that
the sections in the image of $\alpha$ generate ${\cal E}(N)$.
The group $SL(V)$ acts on $\hbox{Hilb}(V \otimes {\cal W}, P)$ and the line
bundle ${\cal L}_m$.
The open subset $Q_2$ is invariant under this action.
\end{say}

Let ${\bf M}_X^\sharp({\cal O}_X, P)$ be the functor for the moduli problem of
$s$-equivalence
classes of $p$-semistable coherent sheaves on $X$ of pure dimension $d$,
Hilbert polynomial $P$, and flat over $S$.

\begin{thm}
\label{thm:simpsonmoduli} \text{{\rm (C. Simpson. \cite{Simpson94})}} $Q_2$ is
contained in
the semistable locus $\Hilb (V \otimes {\cal W}, P)^{ss}({\cal L}_m)$ with
respect to the
action of $SL(V)$ and the linearized line bundle  ${\cal L}_m$ ($m \ge M$).
And the categorical quotient  ${\bf M}_X({\cal O}_X, P)= Q_2/SL(V)$ is a
projective scheme
over $S$ which
coarsely represents the moduli functor  ${\bf M}_X^\sharp({\cal O}_X, P)$.
\end{thm}

We also need  to recall Gieseker's construction of $\overline{M_g}$.

\begin{say}
\label{say:Gieseker's construction} {\sl Gieseker's construction of
$\overline{M_g}$}.
Fix $g \ge 2$, $e= n(2g-2)$ $(n \ge 10)$, $I=e-g$, and a polynomial in $x$,
$p(x) = ex -g +1$.
Set the following Hilbert scheme of subschemes in ${\Bbb P}^I$
$$\hbox{Hilb}_I^{p(x)} :=\{\hbox{subschemes in} \; {\Bbb P}^I \; \hbox{with
Hilbert
polynomial} \; p(x) \}.$$
The group of projective linear transformations $PGL(I+1)$ acts on
$\hbox{Hilb}_I^{p(x)}$
naturally. For the reason of lifting to a linear action, we take $G=SL(I+1)$.
Now consider the locus $H_g$ of n-canonical stable curves in
$\hbox{Hilb}_I^{p(x)}$, that is,
$$H_g =\{[i_{\omega^n}(C)] \in \hbox{Hilb}_I^{p(x)} \}$$
where $C$ is a  DM-stable curve,  $i_{\omega^n}: C \rightarrow {\Bbb P}^I$ is
the
embedding induced by the nth power of the canonical line bundle $\omega$ over
$C$,
and $[i_{\omega^n}(C)]$ is the corresponding Hilbert point of the
n-canonical curve $i_{\omega^n}(C)$. $H_g$ is a $G$-invariant, irreducible,
nonsingular
subscheme of $\hbox{Hilb}_I^{p(x)}$.

By \cite{Gieseker82}, there can be chosen a $SL(I+1)$-linearization on
the Hilbert scheme $\hbox{Hilb}_I^{p(x)}$ such that
\begin{enumerate}
\item  $H_g$ is contained in the stable locus;
\item $H_g$ is closed in the semistable locus; and
\item  the GIT quotient $H_g/SL(I+1)$ is the moduli space $\overline{M_g}$.
\end{enumerate}

For preciseness, we take $e=10(2g-2)$ and $I=10(2g-2) -g$, once and for all.
\end{say}

\begin{say} {\sl The construction of the universal moduli space
${\bf M}_g({\cal O}, P) \rightarrow \overline{M_g}$.}
\label{say:construction}
Let $\widehat{U_g}$ be  the universal curve over
$H_g$.  Consider $X=\widehat{U_g} \rightarrow H_g=S$ as a projective scheme
over the base scheme
$S=H_g$. Fix a Hilbert polynomial $P (x)= r x + d + r(1-g)$.
By Simpson, we get the coarse moduli space   ${\bf M}_X({\cal O}_X, P)$ over
the base scheme
$S=H_g$ of $p$-semistable coherent sheaves (of pure dimension 1)
with the Hilbert polynomial $P (x)= r x + d + r(1-g)$.
 The moduli space   ${\bf M}({\cal O}_X, P)$ (as a projective scheme over
$S=H_g$)
is constructed as the GIT quotient of
$Q_2/S$ by the group $SL(V) =SL(P(N))$ (see \ref{say:Q2} and Theorem
\ref{thm:simpsonmoduli}).
\end{say}

\begin{thm}
\label{thm:universalmoduliofsimpson'sconstructions}
 Fix the Hilbert polynomial $P (x)= r x + d + r(1-g)$.
  The projective
categorical quotient  ${\bf M}_g({\cal O}, P)=(Q_2/H_g)/(SL(V) \times SL(I+1)$
exists
and factors naturally to $\overline{M_g}$ such that over each point $[C] \in
\overline{M_g}$
the fiber is canonically identified with the moduli space  ${\bf M}_C ({\cal
O}_C, P)$ of
$p$-semistable coherent sheaves (of pure dimension $1$) with the Hilbert
polynomial $P$
 modulo the automorphism group of $C$.
\end{thm}
\begin{pf} Consider the map $\pi: Q_2/H_g \rightarrow H_g$ equivariant with
respect to the
projection $\rho: SL(V) \times SL(I+1) \rightarrow SL(I+1)$. For $H_g$ we use a
linearization
$L$ as found by Gieseker (see \ref{say:Gieseker's construction}), while for
$Q_2/H_g$ we use the linearization $\pi^* L^k \otimes {\cal L}_m$ (for $k \gg
0$ and
some sufficiently large $m$, see \ref{say:hilbertschemes}). Now the theorem
follows from
Theorem \ref{thm:generalrgit} or Theorem
\ref{thm:G'toGmainthmforrelativemoduli}.
\end{pf}

\begin{rem}
\label{rem:thesameintheend}
That is, every point of the moduli space ${\bf M}_g({\cal O}, P)$ represents
an equivalence  class of pairs $(C, {\cal E})$ up to automorphism group of $C$,
where
$C$ is a Mumford-Deligne stable curve and ${\cal E}$ is a $p$-semistable
coherent sheaf
of pure dimension $1$ over $C$ with the Hilbert polynomial $P(x) = rx + d+
r(1-g)$.
%A coherent sheaf ${\cal E}$ of rank $r$ is called slope-semistable (resp.
%%stable) if
%$$\frac{\chi({\cal F})
% By the uniqueness of coarse moduli spaces,
%our moduli spaces  ${\bf M}_g({\cal O}, P)$ ought to be isomorphic to the ones
%constructed in \cite{Pandharipande94a}. In particular,  ${\bf M}_g({\cal O},
%%P)$ must
%be isomorphic to the compactified universal Picard $\overline{P^n_g}$ when
% the coherent sheaves ${\cal E}$ are of rank 1
When $P(x) = x + n + 1 - g$ the moduli space ${\bf M}_g({\cal O}, P)$ is
a compactification of the universal Picard $P^n_g$.
It would be interesting to compare our moduli spaces ${\bf M}_g({\cal O}, P)$
with those in \cite{Caporaso94} and \cite{Pandharipande94a}.
%However, as pointed out earlier, our construction
%of the universal moduli spaces  ${\bf M}_g({\cal O}, P)$
%in Theorem \ref{thm:universalmoduliofsimpson'sconstructions}
%bears an easy generalization to the cases when the base moduli spaces are  the
%%moduli spaces
%of higher dimensional varieties (cf. \cite{Kollar90}, \cite{Vieweg95}).
\end{rem}

\section{The universal Hilbert scheme $\hbox{Hilb}^n_g$}
% and the universal Picard $\overline{P^n_g}$}

In this last section, we will give a GIT construction of
the universal Hilbert scheme $\hbox{Hilb}^n_g$
 over $\overline{M_g}$ of $0$-dimensional subschemes of length $n$ on
Mumford-Deligne
stable curves and a canonical morphism from $\hbox{Hilb}^n_g$
to the compactified universal Picard %$\overline{P^n_g} =
$M_g({\cal O}, P)$ where $P(x) = x + n + 1 - g$.

\begin{say} Let $U_g \rightarrow \overline{M_g}$ be the (fake) universal curve
of genus $g \ge 2$
 over $\overline{M_g}$.  $U_g$ has an obvious GIT construction as the quotient
$\widehat{U_g}/SL(I+1)$
by our theorems for $G$-equivariancy RGIT (see \ref{say:construction} for the
definition
of $\widehat{U_g}$). Set $\hbox{Hilb}^n_g \rightarrow \overline{M_g}$
 to be the relative Hilbert scheme over
$\overline{M_g}$ of relatively 0-dimensional subschemes in $U_g$ of length $n$.
Then the fiber of  $\hbox{Hilb}^n_g$ over a point $[C] \in \overline{M_g}$ is
the Hilbert scheme  $\hbox{Hilb}^n_C$ of  0-dimensional subschemes in $C$ of
length $n$ modulo the
automorphism group $\hbox{Aut}(C)$ (we will give a GIT construction of
$\hbox{Hilb}^n_g \rightarrow \overline{M_g}$
in the sequel).
\end{say}

\begin{thm}
\label{thm:mapsfromHilbtopicard}
Let $P(x) = x + n + 1 - g$. Then
there exists a natural dominating morphism $\psi$ from  $\Hilb^n_g$ to
$M_g({\cal O}, P)$ such that
the following diagram is commutative
\begin{equation*}
\begin{CD}
\Hilb^n_g @>{\psi}>> M_g({\cal O}, P) \\
@VVV                         @VVV \\
\overline{M_g} @>{\id}>> \overline{M_g}
\end{CD}
\end{equation*}
\end{thm}

\begin{pf} To construct this morphism scheme-theoretically, we first need to
give a GIT construction
of $\hbox{Hilb}^n_g$ using our theory on $G$-equivariancy RGIT.

Recall that $\overline{M_g}$ is constructed as a GIT quotient of a smooth
irreducible
scheme $H_g$ by the linear transformation $SL(I+1)$ (\ref{say:Gieseker's
construction}).
Let $\widehat{U_g}$ be the  universal family over $H_g$. Set
 $\widehat{\hbox{Hilb}^n_g} \rightarrow H_g$ to be the relative Hilbert scheme
over
$H_g$ of relatively 0-dimensional subschemes in $\widehat{U_g}$ of length $n$.
The group $SL(I+1)$ operates on $\widehat{\hbox{Hilb}^n_g}$ by moving the
subschemes.
Theorems \ref{thm:generalizedReichstein} and
\ref{thm:inducedmapfromG-equivariancy}  imply
 that the GIT quotient $\widehat{\hbox{Hilb}^n_g}/SL(I+1)$ exists
and factors naturally to $H_g/SL(I+1) =  \overline{M_g}$
with the fiber at a point $[C] \in \overline{M_g}$ isomorphic to
$\Hilb^n_C/\hbox{Aut}(C)$. The  quotient
$\widehat{\hbox{Hilb}^n_g}/SL(I+1)$ is our  universal Hilbert scheme
$\hbox{Hilb}^n_g \rightarrow  \overline{M_g}$.

Now given any relatively 0-dimensional subschemes $Z$  in $\widehat{U_g}$ of
length $n$.
We can get a coherent sheave of rank 1 and pure dimension 1,
${\cal O}_{\widehat{U_g}/H_g}(Z) =
{\cal O}_{\widehat{U_g}/H_g}(-Z)^* = (I_Z)^*$.  This leads to a morphism
$$\varphi: \widehat{\hbox{Hilb}^n_g} \longrightarrow Q_2/H_g$$
$$ \varphi: Z \longrightarrow {\cal O}_{\widehat{U_g}/H_g}(Z)$$
(see \ref{say:Q2} for the definition of $Q_2$).
By passing the the quotient we get
$$\varphi': \widehat{\hbox{Hilb}^n_g} \rightarrow Q_2/H_g \rightarrow
(Q_2/H_g)/(SL(V) \times SL(I+1)) = M_g({\cal O}, P)$$
where $P=x + n + 1 -g$.
%(see Theorem \ref{thm:universalmoduliofsimpson'sconstructions}
%and Remark \ref{rem:thesameintheend}).
One checks that $\varphi'$ is constant on the $SL(I+1)$-orbits.
By the universality of categorical quotient, we obtain a canonical morphism
$$\psi:\hbox{Hilb}^n_g =\widehat{\hbox{Hilb}^n_g}/SL(I+1)
\longrightarrow M_g({\cal O}, P).$$  Both $\hbox{Hilb}^n_g$
and $M_g({\cal O}, P)$ are projective and $\psi$ maps surjectively to the
universal
Picard over nonsingular curves of genus $g$. Hnece $\psi$ is dominating.
\end{pf}

\begin{rem} We do not know if the morphism $\psi:\hbox{Hilb}^n_g \rightarrow
M_g({\cal O}, P)$
(or its variants) can be useful in the study of limit linear series for stable
curves
(cf. \cite{EisenbudHarris86}).  We do not know if there exists a canonical
morphism
from the universal FM configuration space over $\overline{M_g}$ to $ M_g({\cal
O}, P)$
so that Theorem \ref{thm:mapsfromHilbtopicard} holds. The
 universal FM configuration space $\overline{FM_{g,n}} \rightarrow
\overline{M_g}$
can be constructed by RGIT as follows. Again let $\widehat{U_g} \rightarrow
H_g$
be the universal family of curves of genus $g$. Let $\widehat{U_g}[n]
\rightarrow H_g $
be the relative  FM configuration space over  $H_g$ (cf.
\cite{Pandharipande94b}).
Then the action of $SL(I+1)$ lifts to a canonical action on $\widehat{U_g}[n]$.
Now applying Theorem \ref{thm:Gmainthmforrelativemoduli} to $\widehat{U_g}[n]
\rightarrow H_g$
and taking the quotients by the group $G=SL(I+1)$, we obtain the
universal FM configuration space
$$\overline{FM_{g,n}} = \widehat{U_g}[n] /SL(I+1) \rightarrow H_g/SL(I+1) =
\overline{M_g}$$
whose fiber at a stable curve $[C]$ is isomorphic to $C[n]/\hbox{Aut}(C)$.
\end{rem}

%\section{The moduli space $\overline{P_g}(r,e,P,\alpha) \rightarrow
%%\overline{M_g}$}
%\label{section:universalmoduli2}

\bibliographystyle{amsplain}
\makeatletter \renewcommand{\@biblabel}[1]{\hfill#1.}\makeatother

\begin{abstract}
We expose in detail the principle that the relative geometric
invariant theory of equivariant  morphisms is related to the
GIT for linearizations near the boundary of the $G$-effective
ample cone. We then apply this principle to construct and
reconstruct various universal moduli spaces. In particular,
we constructed the universal moduli space over $\overline{M_g}$
of Simpson's $p$-semistable coherent sheaves and a canonical
dominating morphism from the universal Hilbert scheme over
$\overline{M_g}$ to a compactified universal Picard.
\end{abstract}

\vskip .4cm

\noindent Department of Mathematics, University of Michigan,  Ann Arbor, MI
48109

\vskip .2cm

\noindent {\sl Electronic mail}:  yihu@@math.lsa.umich.edu

\end{document}